\pdfoutput=1

\documentclass[lettersize,journal]{IEEEtran}
\usepackage{amsmath,amsfonts}
\usepackage{algorithmic}
\usepackage{algorithm}
\usepackage{array}
\usepackage{textcomp}
\usepackage{stfloats}
\usepackage{url}
\usepackage{verbatim}
\usepackage{graphicx}
\usepackage{cite}

\usepackage{color}

\usepackage{amsmath}
\usepackage{amssymb}
\usepackage{mathtools}
\usepackage{amsthm}
\DeclareMathOperator{\tr}{tr}
\DeclareMathOperator{\rank}{rank}
\usepackage{tabularx}
\usepackage{makecell}
\usepackage{booktabs}
\usepackage{caption}
\usepackage{subcaption}
\usepackage{multirow}
\usepackage{calc}  
\usepackage{enumitem}  
\usepackage{pifont}
\usepackage{bm}
\usepackage{relsize}

\newcommand{\cmark}{\ding{51}}%
\newcommand{\xmark}{\ding{55}}%

\newcommand{\Lim}[1]{\raisebox{0.5ex}{\scalebox{0.8}{$\displaystyle \lim_{#1}\;$}}}

\usepackage[capitalize,noabbrev]{cleveref}

\theoremstyle{plain}
\newtheorem{theorem}{Theorem}[section]
\newtheorem{proposition}[theorem]{Proposition}

\newtheorem{corollary}[theorem]{Corollary}
\theoremstyle{definition}
\newtheorem{definition}[theorem]{Definition}
\newtheorem{assumption}[theorem]{Assumption}
\theoremstyle{remark}

\begin{document}
\title{Self-supervised Equality Embedded Deep Lagrange Dual for Approximate Constrained Optimization}
\author{Minsoo Kim, Hongseok Kim,~\IEEEmembership{Senior member,~IEEE,}
\thanks{This paper has been published in IEEE Transactions on Power Systems \cite{kim2024unsupervised}. The published version includes simulation results on a 1354-bus system and provides a more detailed analysis from a power systems perspective.}
}


\maketitle
\begin{abstract}
Conventional solvers are often computationally expensive for constrained optimization, particularly in large-scale and time-critical problems. While this leads to a growing interest in using neural networks (NNs) as fast optimal solution approximators, incorporating the constraints with NNs is challenging. In this regard, we propose deep Lagrange dual with equality embedding (\texttt{DeepLDE}), a framework that learns to find an optimal solution without using labels. To ensure feasible solutions, we embed equality constraints into the NNs and train the NNs using the primal-dual method to impose inequality constraints. Furthermore, we prove the convergence of \texttt{DeepLDE} and show that the primal-dual learning method alone cannot ensure equality constraints without the help of equality embedding. Simulation results on convex, non-convex, and AC optimal power flow (AC-OPF) problems show that the proposed \texttt{DeepLDE} achieves the smallest optimality gap among all the NN-based approaches while always ensuring feasible solutions. Furthermore, the computation time of the proposed method is about 5 to 250 times faster than \texttt{DC3} and the conventional solvers in solving constrained convex, non-convex optimization, and/or AC-OPF.
\end{abstract}


\IEEEpeerreviewmaketitle

\section*{Nomenclature}
\addcontentsline{toc}{section}{Nomenclature}

\begin{IEEEdescription}[\IEEEusemathlabelsep\IEEEsetlabelwidth{$\hat{\mathcal{S}}_k$}]
\item[$\mathcal{L}$] Lagrange function
\item[$\mathcal{L}_e$] Lagrange function with equality embedding
\item[$\mathcal{D}$] Training dataset
\item[$\mathbf{d}$] Input data
\item[$f_\mathbf{d}$] Objective function
\item[$g_\mathbf{d}$] Inequality constraints
\item[$h_\mathbf{d}$] Equality constraints
\item[$\bm{\lambda}$] Lagrange multiplier of $g_\mathbf{d}$
\item[$\bm{\mu}$] Lagrange multiplier of $h_\mathbf{d}$
\item[$\Phi$] Neural network model
\item[$\mathbf{w}$] Model parameters of $\Phi$  
\item[$\Psi_h$] Equality embedding function
\item[$\mathbf{y}$] Output of $\Phi(\mathbf{w},\mathbf{d})$ without using equality embedding
\item[$\mathbf{x}$] Output of $\Phi(\mathbf{w},\mathbf{d})$ when using equality embedding
\item[$\mathbf{z}$] Output of $\Psi_h(\mathbf{x})$
\end{IEEEdescription}

\section{Introduction}
\IEEEPARstart{S}{olving} constrained optimization with conventional solvers are often expensive. The iterative methods, widely used in conventional solvers, are usually built on concrete mathematical theories about convergence \cite{boyd2004convex}, but computation time is critical in time-sensitive large-scale constrained optimization. Hence, using neural networks (NNs) as fast approximators for optimization problems to obtain feasible and optimal solutions is of great interest \cite{hasan2020survey, bengio2021machine}.

While there is significant potential in learning constrained optimization with NNs \cite{bengio2021machine}, several challenges exist. \textbf{First}, imposing constraints on NNs and obtaining a feasible solution is non-trivial (feasibility). \textbf{Second}, it is desired that the NNs learn to solve constrained optimization without using ready-made solutions as labels (self-supervised). \textbf{Third}, NNs should approximate constrained optimization significantly faster than conventional solvers (computation time). 

There have been several efforts in the literature to resolve these three challenges. As regards the first challenge, the penalty terms of constraint violations were integrated into the loss functions in \cite{chen2022learning, pan2020deepopf, nellikkath2021physics, nellikkath2022physics}. To reduce the burden of manually tuning the weights of the penalty terms, the primal-dual method was applied to impose constraints to NNs in \cite{nandwani2019primal, fioretto2020Lagrangian, fioretto2020predicting, kotary2022fast, park2022self}. However, some of these works did not overcome the second challenge since they used pre-solved solutions as labels \cite{nandwani2019primal, pan2020deepopf, fioretto2020Lagrangian, fioretto2020predicting, nellikkath2021physics,nellikkath2022physics}. In addition, the penalty and the primal-dual method reveal limited capabilities in ensuring equality constraints.

Instead of relying on the penalty or the primal-dual method to ensure feasibility, the works of \cite{chen2021enforcing, kim2022projection} embedded differentiable solvers \cite{amos2017optnet, agrawal2019differentiable} into NNs. However, their approaches are limited to solve only specific types of constrained optimization. Furthermore, employing differentiable solvers in the forward path of NNs in training and inference incurred long computation time and confronted the third challenge.

To mitigate these limitations, \cite{donti2021DC3} proposed a deep constraint completion and correction (\texttt{DC3}) method. Using the completion and the correction methods, \texttt{DC3} achieved the feasibility of the inferred solutions and fast computation without using additional pre-solved solutions. However, the correction method needs multiple iterations, which incurs excessive computation time. Furthermore, inefficient memory usage due to the Jacobian matrix of inequality constraints for completion is another drawback of \texttt{DC3}.

In this regard, we propose a novel technique called deep Lagrange dual with equality embedding (\texttt{DeepLDE}) to overcome the aforementioned three challenges all together. The NN \textit{part} of \texttt{DeepLDE} predicts a part of optimization solution. Then, the rest of the solution is determined by embedded linear or nonlinear equality constraints. We employ implicit differentiation to backpropagate through the embedded equality constraints. To satisfy inequality constraints, we apply the primal-dual method which substantially reduces the computation time. 

We summarize our key contributions as follows:\IEEEpubidadjcol

First, we provide a novel self-supervised learning framework, called \texttt{DeepLDE}, that approximates non-convex optimization solutions using an equality embedding method in conjunction with the Lagrange dual method. To the best of our knowledge, we first demonstrate that it works in self-supervised learning. To ensure both equality and inequality constraints, we train the NN in an alternating manner: minimizing the loss function with equality embedding in the inner loop and maximizing the dual function in the outer loop. Interestingly, however, this nested loop of training does not increase the number of epochs (see the learning curves in Fig.~\ref{fig:EE_comparison}) and ensures that the solution provided by \texttt{DeepLDE} remains feasible.

Second, \texttt{DeepLDE} reduces the inference time by 80\% compared to \texttt{DC3}, the state-of-the-art self-supervised non-convex optimization solver. This is achieved by eliminating the time-consuming correction procedure used by \texttt{DC3} in satisfying inequality constraints. We also compare \texttt{DeepLDE} with its supervised learning version (\texttt{SL+LDE}) and find that the proposed self-supervised learning performs as well as or \textit{even better} than its supervised counterpart, without the need for pre-labeled training data.

Third, by leveraging the concept of $\varepsilon$-inexact solution (see Definition~\ref{def:eps_inexact}), we show that equality constraints cannot be met without the use of equality embedding (see Proposition~\ref{prop:inexact_equality} and Corollary~\ref{cor:positive_error}).

Fourth, we rigorously verify \texttt{DeepLDE} against ten state-of-the-art algorithms including \texttt{OSQP}, \texttt{qpth}, and \texttt{IPOPT} in solving both convex and non-convex optimization problems. Furthermore, using IEEE standard 57 and 118 bus systems \cite{babaeinejadsarookolaee2019power}, we extensively demonstrate that \texttt{DeepLDE} outperforms \texttt{PYPOWER} and \texttt{DC3}, which are designed to solve large-scale real power system problems, in terms of feasibility, optimality, and/or computation time.

It should be noted that the proposed method is in part related to the works of \cite{fioretto2020Lagrangian, fioretto2020predicting, donti2021DC3} but differs in several key aspects. First, in contrast to \cite{fioretto2020Lagrangian, fioretto2020predicting}, \texttt{DeepLDE} employs the equality embedding method and achieves zero-constraint violations. Second, instead of using the time-consuming correction method of \cite{donti2021DC3}, our work leverages the primal-dual method to expedite imposing inequality constraints on NNs. Unlike the correction in \cite{donti2021DC3}, which is applied as a post-processing during inference, the primal-dual method in \texttt{DeepLDE} is used solely in the training process.

The rest of the paper is structured as follows. Section~\ref{sec:Related_work} provides an overview of related works while Section~\ref{sec:Preliminaries} presents the problem formulation. The methodology in \texttt{DeepLDE} and mathematical analysis are described in Section~\ref{sec:Methodology}, and the simulation settings and results are discussed in Section~\ref{sec:experi}. Finally, Section~\ref{sec:conclusion} concludes with a discussion on future research directions.

\section{Related Work}
\label{sec:Related_work}
Early works of learning in optimization focused on hyperparameter optimization \cite{bergstra2011algorithms, hutter2011sequential, snoek2012practical, bergstra2012random}. Recently, interactions between learning and optimization have been tried in diverse fields in machine learning. For instance, the works of \cite{li2017learning, chenL2O2022learning} reduced iterations of conventional solvers by learning an optimization algorithm from the given dataset or task. The domain-specific learning method was proposed by \cite{donti2017task,wilder2019end,poganvcic2019differentiation,elmachtoub2022smart}. Since the computation of trained NNs is faster than that of conventional solvers \cite{pan2020deepopf}, there have been increasing efforts to integrate constraints with NNs \cite{kotary2021end}.

\textbf{Using activation functions for simple constraints.} Incorporating simple constraints through the use of activation functions is a common practice in learning constrained optimization. For example, using ReLU or Sigmoid after the last layer ensures that the outputs of NNs are always non-negative \cite{sun2018learning} or bounded \cite{zamzam2020learning}. However, these approaches have limited ability to incorporate equality and more intricate inequality constraints.

\textbf{Penalty method.} 
Several works integrated penalty terms into the loss function to impose constraints with the NNs. The works of \cite{pan2020deepopf} defined the penalty terms as $l_1$- or $l_2$-norm of the constraint violation degree. The works of \cite{zhang2021convex, nellikkath2021physics, chen2022learning, nellikkath2022physics} used the KKT conditions as penalty terms to obtain optimal and feasible solutions from NNs. Although the penalty methods reduce the constraint violations, the feasibility of the solutions is not ensured. \cite{pan2020deepopf} projected the outputs of a trained NN onto the feasible region, but this method is limited to convex feasible sets. Furthermore, \cite{venzke2020learning, ul2022learning} examined the worst-case constraint violations when imposing constraints on neural networks, but they did not minimize the worst-case violations. Lastly, \cite{zhao2023ensuring} proposed preventive learning to ensure the solution of NNs remains feasible, but this approach was limited to linear constraints only.

\textbf{Primal-dual learning method.}
One of the major drawbacks of the penalty method is that manually selecting the weights of penalty terms is time-consuming. Thus, several works applied the primal-dual learning method \cite{nandwani2019primal, fioretto2020Lagrangian, fioretto2020predicting} and reduced the burden of manual tuning. The work of \cite{kotary2022fast} solved a combinatorial optimization problem using the primal-dual learning method. The work in \cite{park2022self} used two NNs, i.e., one for determining the primal variables and the other for the dual variables. However, the primal-dual learning method cannot solely ensure equality constraints.

\textbf{Non-penalty method.} 
Instead of incorporating constraints into the loss function, the authors of \cite{ling2020solving, ling2021can, tabas2022computationally, li2022learning} directly mapped the inputs of the NN onto a feasible region. However, there were some limitations to these approaches. For example, \cite{ling2020solving, ling2021can} required pre-solved labels for training, and \cite{tabas2022computationally, li2022learning} needed a specific feasible region \cite{blanchini2008set} and at least one interior point. Specifically, the work in \cite{li2022learning} obtained the interior point by solving linear programming, but their approach showed limitations in solving large-scale constrained optimization problems. Additionally, \cite{he2015neural} employed the barrier Lyapunov function and NN to handle output constraints, but their approach was limited to inequality constraints.

\textbf{Implicit layers.}
Recent research has focused on embedding implicit layers within NNs in order to construct implicit functions between inputs and outputs of NNs \cite{amos2017optnet,agrawal2019differentiable,donti2021DC3}. In particular, the implicit layers were employed to embed optimization problems in NNs. For example, the works in \cite{amos2017optnet} integrated a QP solver into the NN, and the authors of \cite{agrawal2019differentiable} generalized the works of \cite {amos2017optnet} by embedding convex optimization layers. However, the computation time of their approaches was significantly longer than naive NNs. \texttt{DC3} \cite{donti2021DC3} reduced the computation time but still suffered from the additional computational overhead of the correction method to satisfy the inequality constraints.

\section{Problem Formulation}
\label{sec:Preliminaries}
\begin{figure}[t]
	\centering
	\includegraphics[width=0.99\columnwidth]{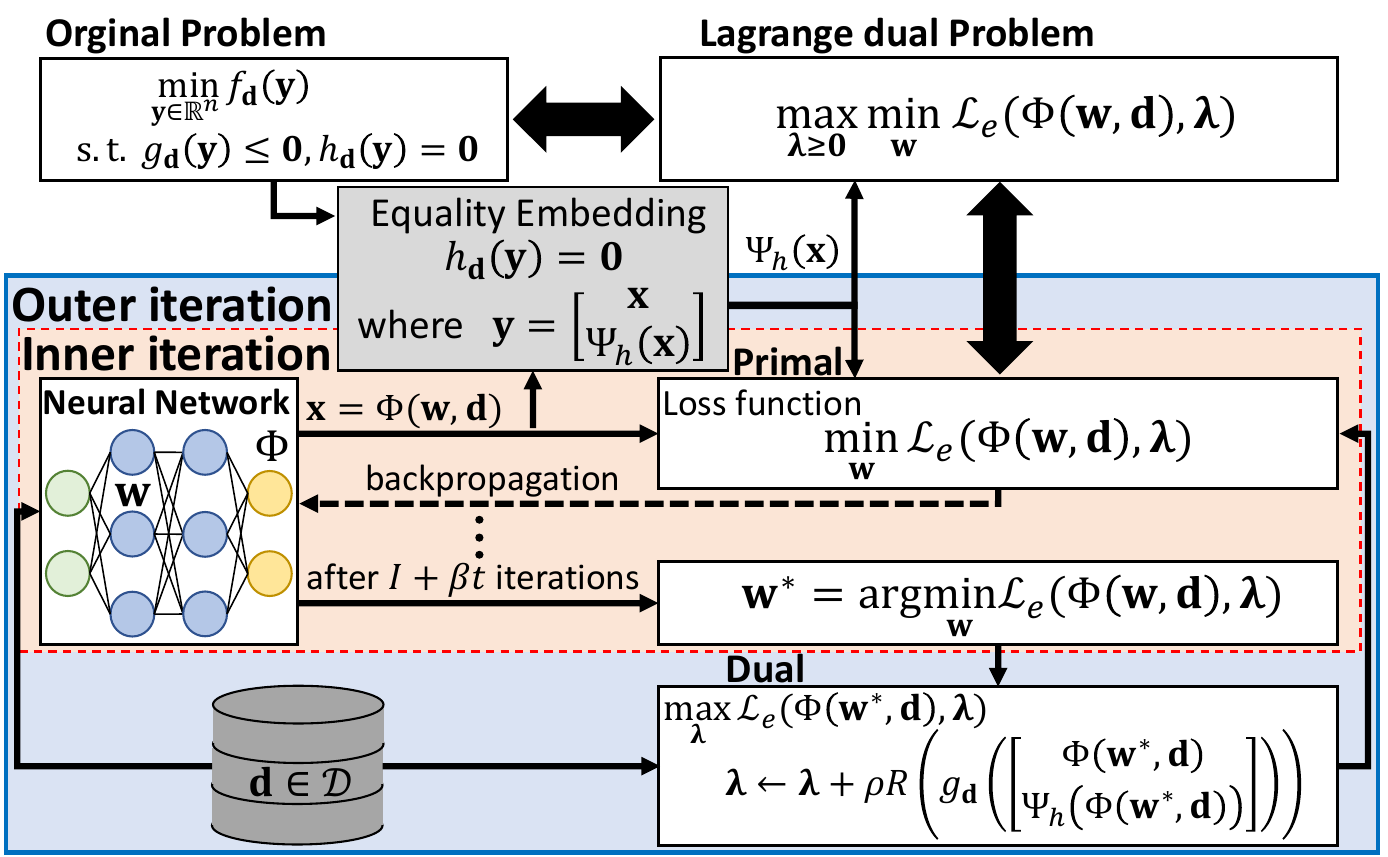}
	\caption{Overall framework of \texttt{DeepLDE}.}
	\label{fig:DeepLDE_framework}

\end{figure}
We consider families of constrained optimization problems; given each input data point $\mathbf{d}\in\mathbb{R}^n$ in $\mathcal{D}$, predict (or learn) an optimal solution $\mathbf{y}$ of solving
\begin{equation}
\begin{aligned}
    &\min_{\mathbf{y}\in\mathbb{R}^n} f_{\mathbf{d}}(\mathbf{y}),\\
    &\textcolor{black}{\text{s.t.}} \quad g_{\mathbf{d}}(\mathbf{y})\leq \mathbf{0},\quad h_{\mathbf{d}}(\mathbf{y}) = \mathbf{0},
    \label{eq:prelim_opt}
\end{aligned}
\end{equation}
where the objective function $f_\mathbf{d}:\mathbb{R}^{n}\rightarrow\mathbb{R}$, the inequality constraints $g_{\mathbf{d}}:\mathbb{R}^n\rightarrow\mathbb{R}^{n_{\text{ineq}}}$ and the equality constraints $h_{\mathbf{d}}:\mathbb{R}^n\rightarrow\mathbb{R}^{n_{\text{eq}}}$ are potentially non-convex.

In solving (\ref{eq:prelim_opt}), let $\bm{\lambda} \in \mathbb{R}^{n_{\text{ineq}}}$ and $\bm{\mu}\in\mathbb{R}^{n_{\text{eq}}}$ denote the Lagrange multiplier of $g_{\mathbf{d}}$ and $h_{\mathbf{d}}$, respectively, where $\bm{\lambda}\geq0$. Then, we have the following Lagrangian function:
\begin{equation}
    \mathcal{L}(\mathbf{y}, \bm{\lambda}, \bm{\mu}) = f_{\mathbf{d}}(\mathbf{y}) + \bm{\lambda}^{\intercal}R(g_{\mathbf{d}}(\mathbf{y})) + \bm{\mu}^{\intercal}S(h_{\mathbf{d}}(\mathbf{y})),
    \label{eq:prelim_lag}
\end{equation}
where $R(g_{\mathbf{d}}(\cdot)) = \max\{g_{\mathbf{d}}(\cdot),0\}$, and $S(h_{\mathbf{d}}(\cdot)) = |h_{\mathbf{d}}(\cdot)|$. Here, $\max\{\cdot, 0\}$ and $|\cdot|$ are elementwise operators that return non-negative and absolute values, respectively. Using the Lagrangian function in (\ref{eq:prelim_lag}), the optimization problem (\ref{eq:prelim_opt}) becomes
\begin{equation}
    \max_{\bm{\mu}}\max_{\bm{\lambda}\geq 0}\min_{\mathbf{y}}\mathcal{L}(\mathbf{y}, \bm{\lambda}, \bm{\mu}).
    \label{eq:prelim_maxmin}
\end{equation}
For a given input vector $\mathbf{d}$, we solve (\ref{eq:prelim_maxmin}) using a subgradient method that updates $\mathbf{y}$, $\bm{\lambda}$, and $\bm{\mu}$ iteratively as follows \cite{fontaine2014constraint}:
\begin{equation}
    \mathbf{y}^* = \arg\min_{\mathbf{y}}\mathcal{L}(\mathbf{y}, \bm{\lambda}, \bm{\mu}),
    \label{eq:prelim_primal}
\end{equation}
\begin{equation}
    \bm{\lambda} \leftarrow \bm{\lambda} + \rho\nabla_{\bm{\lambda}} \mathcal{L}(\mathbf{y}^*, \bm{\lambda}, \bm{\mu}),
    \label{eq:prelim_dual_ineq}
\end{equation}
\begin{equation}
    \bm{\mu} \leftarrow \bm{\mu} + s\nabla_{\bm{\mu}} \mathcal{L}(\mathbf{y}^*, \bm{\lambda}, \bm{\mu}),
    \label{eq:prelim_dual_eq}
\end{equation}
where $\rho>0$ and $s>0$ are the step sizes of $\bm{\lambda}$ and $\bm{\mu}$, respectively.

\section{Methodology}
\label{sec:Methodology}
\subsection{Embedding Equality Constraints}
Instead of using the Lagrange multiplier $\bm{\mu}$ for $S(h_{\mathbf{d}}(\cdot))$, we embed the equality constraints in the original problem \cite{nocedal2006numerical}. We assume that $h_\mathbf{d}:\mathbb{R}^n\rightarrow \mathbb{R}^{n_{\text{eq}}}$ is not overdetermined, i.e., $n\geq n_{\text{eq}}$. Thus, if $n-n_{\text{eq}}$ elements of $\mathbf{y}$ are computed, the rest of $n_{\text{eq}}$ elements that satisfy $h_{\mathbf{d}}(\mathbf{y}) = \mathbf{0}$ can be accordingly determined. Thus, we decompose $\mathbf{y}=[\mathbf{x}^{\intercal}\;\;\mathbf{z}^\intercal]^{\intercal}$ such as $\mathbf{x}\in\mathbb{R}^{n-n_{\text{eq}}}$ and $\mathbf{z}\in\mathbb{R}^{n_{\text{eq}}}$ satisfying
\begin{equation}
    h_{\mathbf{d}}(\mathbf{y}) = 
        h_{\mathbf{d}}\Big(\begin{bmatrix}
        \mathbf{x}\\
        \mathbf{z}
    \end{bmatrix}\Big) =
        h_{\mathbf{d}}\Big(\begin{bmatrix}
        \mathbf{x}\\
        \Psi_{h}(\mathbf{x})
    \end{bmatrix}\Big) = \mathbf{0},
    \label{eq:method_separateY}
\end{equation}
where $\Psi_h:\mathbb{R}^{n-n_{\text{eq}}}\rightarrow\mathbb{R}^{n_{\text{eq}}}$ is a function that satisfies the equality constraints given $\mathbf{x}$, i.e., $\mathbf{z} = \Psi_h(\mathbf{x})$. Now, the optimization problem (\ref{eq:prelim_opt}) is reformulated without the equality constraints:
\begin{equation}
\begin{aligned}
    &\min_{\mathbf{x}\in\mathbb{R}^{n-n_{\text{eq}}}} f_{\mathbf{d}}\Big(
    \begin{bmatrix}
        \mathbf{x}\\
        \Psi_{h}(\mathbf{x})
    \end{bmatrix}\Big),\\
    &\textcolor{black}{\text{s.t.}} \quad g_{\mathbf{d}}\Big(
    \begin{bmatrix}
        \mathbf{x}\\
        \Psi_{h}(\mathbf{x})
    \end{bmatrix}\Big)\leq 0.
    \label{eq:method_reformulZ}
\end{aligned}
\end{equation}
Then, the equality embedded Lagrangian function of (\ref{eq:method_reformulZ}) is
\begin{equation}
\begin{aligned}
    \mathcal{L}_e(\mathbf{x}, \bm{\lambda})
    &= f_{\mathbf{d}}\Big(\begin{bmatrix}
        \mathbf{x}\\
        \Psi_{h}(\mathbf{x})
    \end{bmatrix}\Big)+\bm{\lambda}^{\intercal}R\bigg(g_{\mathbf{d}}\Big(\begin{bmatrix}
        \mathbf{x}\\
        \Psi_{h}(\mathbf{x})
    \end{bmatrix}\Big)\bigg).
    \label{eq:method_lagElim}
\end{aligned}
\end{equation}
\begin{algorithm}[t]
\caption{Training of \texttt{{DeepLDE}}}\label{alg:DeepLDE_noW}
\begin{algorithmic}[1]
\STATE \textbf{Inputs :}\\Maximum number of outer iterations $T$, Maximum number of inner iterations $I$, Training dataset $\mathcal{D}$, Learning rate $\eta$, Step size $\rho_0$, Initial value of the Lagrange multiplier $\bm{\lambda}_0$, Increments of iteration $\beta$, Decrements of step size $\gamma$
\STATE \textbf{Initialize :}\\Randomly initialize $\mathbf{w}$
\STATE $(\rho, \bm{\lambda}) \gets (\rho_0, \bm{\lambda}_0)$
\FOR{$t \in \{1, ..., T\}$}
\FOR{$i \in \{1, ..., I\}$}
\FORALL{$\mathbf{d} \in\mathcal{D}$}
\STATE $\mathbf{x} \gets \Phi(\mathbf{w}, \mathbf{d})$
\STATE $\nabla_{\mathbf{w}}\mathcal{L}_e(\mathbf{x}, \bm{\lambda}) \gets \text{calculate via (\ref{eq:gradient_elim}) using }\mathbf{x}$
\STATE $\mathbf{w} \gets \mathbf{w} - \eta\nabla_\mathbf{w}\mathcal{L}_e(\mathbf{x}, \bm{\lambda})$
\ENDFOR
\ENDFOR
\STATE $\bm{\lambda} \gets \bm{\lambda} + \rho\mathlarger{\sum}_{\mathbf{d}\in\mathcal{D}}R\bigg(g_{\mathbf{d}}\Big(
        \begin{bmatrix}
        \mathbf{\Phi({\mathbf{w}}, \mathbf{d})}\\
        \Psi_h(\mathbf{\Phi(\mathbf{w}, \mathbf{d}))}
    \end{bmatrix}\Big)\bigg)$
\STATE $I \gets I + \beta$
\STATE $\rho \gets \dfrac{\rho_0}{1+\gamma t}$
\ENDFOR
\end{algorithmic}
\end{algorithm}

\subsection{Learning Constrained Optimization with NNs}
Let $\Phi$ denote an NN model. Instead of directly predicting the entire elements of $\mathbf{y}$, we train $\Phi$ to predict an optimal $\mathbf{x}$. After that, we determine $\mathbf{y}$ using $\Psi_h$ as shown in (\ref{eq:method_separateY}). Let $\mathbf{w}$ be the model parameters of $\Phi$. Then, given $\mathbf{d}$ as an input of $\Phi$, we have
\begin{equation}
    \mathbf{x} = \Phi(\mathbf{w}, \mathbf{d}).
    \label{eq:defineZ_NN}
\end{equation}
 Now, we construct the Lagrangian function by plugging in (\ref{eq:defineZ_NN}) into (\ref{eq:method_lagElim}):
\begin{equation}
\begin{aligned}
    \mathcal{L}_e(\Phi(\mathbf{w}, \mathbf{d}), \bm{\lambda}) &= f_\mathbf{d}\Big(
    \begin{bmatrix}
        \Phi(\mathbf{w}, \mathbf{d})\\
        \Psi_h(\Phi(\mathbf{w}, \mathbf{d}))
    \end{bmatrix}\Big)\\
    &+ \bm{\lambda}^{\intercal}R\bigg(g_\mathbf{d}\Big(
    \begin{bmatrix}
        \Phi(\mathbf{w}, \mathbf{d})\\
        \Psi_h(\Phi(\mathbf{w}, \mathbf{d}))
    \end{bmatrix}\Big)\bigg),
\end{aligned}
\label{eq:method_LagW}
\end{equation}
which is optimized by the following max-min optimization:
\begin{equation}
    \max_{\bm{\lambda}\geq 0}\min_{\mathbf{w}}\mathcal{L}_e.
    \label{eq:method_maxminW}
\end{equation}
We again employ the subgradient method to solve (\ref{eq:method_maxminW}) \cite{fontaine2014constraint} and update $\mathbf{w}$ and $\bm{\lambda}$ alternatively,
\begin{equation}
    \mathbf{w}^* = \arg\min_{\mathbf{w}}\mathcal{L}_e,
    \label{eq:updateW}
\end{equation}
\begin{equation}
    \bm{\lambda} \leftarrow \bm{\lambda} + \rho\nabla_{\bm{\lambda}}\mathcal{L}_e,
    \label{eq:updateLambda}
\end{equation}
where $\rho>0$ is the step size of $\bm{\lambda}$, and $\nabla_{\bm{\lambda}}\mathcal{L}_e$ is given by
\begin{equation}
    \nabla_{\bm{\lambda}}\mathcal{L}_e = R\bigg(g_\mathbf{d}\Big(
    \begin{bmatrix}
        \Phi(\mathbf{w}^*, \mathbf{d})\\
        \Psi_h(\Phi(\mathbf{w}^*, \mathbf{d}))
    \end{bmatrix}\Big)\bigg).
\end{equation}
The approximation of (\ref{eq:updateW}) is done by the stochastic gradient descent (SGD) method as done in \cite{fioretto2020Lagrangian,fioretto2020predicting,kotary2022fast,park2022self}.
If the embedded equality constraints are nonlinear, $\Psi_h$ denotes a procedure (e.g., Newton's method) to solve the nonlinear equation, and implicit differentiation is used to update the parameters of NN by backpropagating the loss function through $\Psi_h$ \cite{amos2017optnet,donti2021DC3}. In doing this, the chain rule lets us have
\begin{equation}
    \dfrac{d}{d\mathbf{x}}h_{\mathbf{d}}\Big(\begin{bmatrix}
        \mathbf{x}\\
        \Psi_h(\mathbf{x})
    \end{bmatrix}\Big) = \dfrac{\partial h_{\mathbf{d}}}{\partial\mathbf{x}} + \dfrac{\partial h_{\mathbf{d}}}{\partial\Psi_h(\mathbf{x})}\dfrac{{\partial\Psi_h(\mathbf{x})}}{\partial \mathbf{x}} = 0,
    \label{eq:method_implicit}
\end{equation}
which implies
\begin{equation}
    \dfrac{{\partial\Psi_h(\mathbf{x})}}{\partial \mathbf{x}} = -\Big(\dfrac{\partial h_{\mathbf{d}}}{\partial\Psi_h(\mathbf{x})}\Big)^{-1}\dfrac{\partial h_{\mathbf{d}}}{\partial\mathbf{x}}.
    \label{eq:method_implicit2}
\end{equation}
Thus, the backpropagation process to update $\mathbf{w}$ is
\begin{equation}
\begin{aligned}
    \nabla_{\mathbf{w}}\mathcal{L}_e &= \Big(\dfrac{\partial\mathcal{L}_e}{\partial\mathbf{x}} + \dfrac{\partial\mathcal{L}_e}{\partial\Psi_h(\mathbf{x})}\dfrac{\partial\Psi_h(\mathbf{\mathbf{x}})}{\partial\mathbf{x}}\Big)\dfrac{\partial\mathbf{x}}{\partial\mathbf{w}}\\&=
    \Big(\dfrac{\partial\mathcal{L}_e}{\partial\mathbf{x}} - \dfrac{\partial\mathcal{L}_e}{\partial\Psi_h(\mathbf{x})}\Big(\dfrac{\partial h_{\mathbf{d}}}{\partial\Psi_h(\mathbf{x})}\Big)^{-1}\dfrac{\partial h_{\mathbf{d}}}{\partial\mathbf{x}}\Big)\dfrac{\partial\mathbf{x}}{\partial\mathbf{w}}.
\end{aligned}
    \label{eq:gradient_elim}
\end{equation}

We summarize the overall \texttt{DeepLDE} framework in Fig.~\ref{fig:DeepLDE_framework} and the training process in Algorithm~ \ref{alg:DeepLDE_noW}. Next, we provide the following proposition to ensure the convergence of Algorithm~\ref{alg:DeepLDE_noW}:
\begin{proposition}
Let $\beta\in\mathbb{Z}$ and $\gamma\in\mathbb{R}$ be nonnegative values. Then, Algorithm~\ref{alg:DeepLDE_noW} converges if $\beta > 0$ or $\gamma > 0$.
\label{prop:local_convergence}
\end{proposition}
The proof of Proposition~\ref{prop:local_convergence} is in Appendix~\ref{subsec:proof_converge}. The intuition of Proposition~\ref{prop:local_convergence} is to approximate (\ref{eq:updateW}) to update $\mathbf{w}$ for a small perturbation of $\bm{\lambda}$ \cite{jin2020local}. The convergence analysis when $\beta>0$ is provided in \cite{nandwani2019primal}, i.e., when the number of iterations increases in epochs. We add an additional sufficient condition for convergence, i.e., $\gamma>0$; the step size of $\rho$ diminishes in epochs.

\subsection{Analysis about Equality Constraint Embedding}
\begin{algorithm}[t]
\caption{\texttt{DeepLDE} without equality embedding (denoted by \texttt{LDF} hereafter)}\label{alg:LDF_noW}
\begin{algorithmic}[1]
\STATE \textbf{Inputs :}\\Maximum number of outer iterations $T$, Maximum number of inner iterations $I$, $\mathbf{d}$ in training dataset $\mathcal{D}$, Learning rate $\eta$, Step size $\rho_0$ and $s_0$, Initial value of the Lagrange multiplier $\bm{\lambda}_0$ and $\bm{\mu}_0$, Increments of outer iteration $\beta$, Decrements of step size $\gamma$
\STATE \textbf{Initialize :}\\Randomly initialize $\mathbf{w}$
\STATE $(\rho, s, \bm{\lambda}, \bm{\mu}) \gets (\rho_0, s, \bm{\lambda}_0, \bm{\mu}_0)$
\FOR{$t \in \{1, ..., T\}$}
\FOR{$i \in \{1, ..., I\}$}
\FORALL{$\mathbf{d} \in\mathcal{D}$}
\STATE $\mathbf{y} \gets \Phi(\mathbf{w}, \mathbf{d})$
\STATE $\mathbf{w} \gets \mathbf{w} - \eta\nabla_{\mathbf{w}}\mathcal{L}(\mathbf{y}, \bm{\lambda}, \bm{\mu})$
\ENDFOR
\ENDFOR
\STATE $\bm{\lambda} \gets \bm{\lambda} + \rho
\mathlarger{\sum}_{\mathbf{d}\in\mathcal{D}}R(g_{\mathbf{d}}(\Phi(\mathbf{w}, \mathbf{d})))$
\STATE $\bm{\mu} \gets \bm{\mu} + s
\mathlarger{\sum}_{\mathbf{d}\in\mathcal{D}}S(h_{\mathbf{d}}(\Phi(\mathbf{w}, \mathbf{d})))$
\STATE $I \gets I + \beta$
\STATE $\rho \gets \dfrac{\rho_0}{1+\gamma t}$
\STATE $s \gets \dfrac{s_0}{1+\gamma t}$
\ENDFOR
\end{algorithmic}
\end{algorithm}
Now, we analyze the rationale and benefits of using equality embedding in the proposed method. In doing this, we provide Algorithm~\ref{alg:LDF_noW}, which replaces equality embedding within \texttt{DeepLDE} by the update of $\bm{\mu}$ in (\ref{eq:prelim_dual_eq}) for analytical comparison. So, the output vector of NN is $\mathbf{y}\in\mathbb{R}^n$ as done in (\ref{eq:prelim_primal})$-$(\ref{eq:prelim_dual_eq}) where (\ref{eq:prelim_primal}) is approximated by SGD. Note that Algorithm~\ref{alg:LDF_noW} is similar to \cite{fioretto2020Lagrangian, fioretto2020predicting} but does not use labels. This method is denoted by Lagrange dual framework (\texttt{LDF}) hereafter.
\begin{definition}[$\varepsilon$-inexact solution]
Let $\mathbf{y}^*$ be an optimal solution of (\ref{eq:prelim_opt}). Then, $\mathbf{y}_{\varepsilon}$, the solution predicted from Algorithm~\ref{alg:LDF_noW}, is said $\varepsilon$-inexact when $\mathbf{y}_{\varepsilon} = \mathbf{y}^* + \varepsilon$.
\label{def:eps_inexact}
\end{definition}

\begin{assumption} (normally distributed errors). $\varepsilon$-inexact solution is normally distributed, i.e., $\varepsilon\sim\mathcal{N}(\mathbf{0}, \delta^2\mathbf{I})$.
\label{assume:normal_error}
\end{assumption}
Using Definition~\ref{def:eps_inexact} and Assumption~\ref{assume:normal_error}, we provide Proposition~\ref{prop:inexact_equality} and Corollary~\ref{cor:positive_error}, which imply that $h_{\mathbf{d}}(\mathbf{y}_\varepsilon)=\mathbf{0}$ in (\ref{eq:prelim_opt}) cannot hold when Algorithm~\ref{alg:LDF_noW} is used. For notational simplicity, we use $h = h_{\mathbf{d}}$ hereafter.
\begin{proposition}
Let $h:\mathbb{R}^{n}\rightarrow\mathbb{R}^{n_{\text{eq}}}$ be a continuously differentiable function, and $\mathbf{J}_h(\mathbf{y})$ be the Jacobian matrix of $h$. Then, for any small $\varepsilon$, we have
\begin{equation}
    \mathbb{E}[\mathbf{1}^{\intercal}S(h(\mathbf{y}_\varepsilon))] = \sqrt{\dfrac{2}{\pi}}\sigma||\mathbf{J}_h(\mathbf{y}^*)||_*
    \label{eq:prop_inequality}
\end{equation}
where $||\mathbf{J}_h||_* = \tr(\sqrt{\mathbf{J}_h\mathbf{J}_h^{\intercal}})$.
\label{prop:inexact_equality}
\end{proposition}
The proof of Proposition~\ref{prop:inexact_equality} is provided in Appendix~\ref{subsec:inexact_equality}. Proposition~\ref{prop:inexact_equality} says that $\mathbb{E}[\mathbf{1}^{\intercal}S(h(\mathbf{y}_\varepsilon))]$ is determined by $\sigma$ and $||\mathbf{J}_h||_*$. From Proposition~\ref{prop:inexact_equality}, we have Corollary~\ref{cor:positive_error}.
\begin{corollary}
    Let $h:\mathbb{R}^{n}\rightarrow\mathbb{R}^{n_{\text{eq}}}$ be a continuously differentiable function, and $\mathbf{J}_h(\mathbf{y})$ be the Jacobian matrix of $h$. Then, for any small $\varepsilon$, if $\mathbf{J}_h(\mathbf{y}^*) \neq \mathbf{0}$, we have
    \begin{equation}
        \mathbb{E}[\mathbf{1}^{\intercal}S(h(\mathbf{y}_\varepsilon))]>0.
    \end{equation}
    \label{cor:positive_error}
\end{corollary}

The proof of Corollary~\ref{cor:positive_error} is provided in Appendix~\ref{subsec:positive_error}. The condition $\mathbf{J}_h(\mathbf{y}^*) \neq \mathbf{0}$ in Corollary~\ref{cor:positive_error} often occurs because $h(\mathbf{y}^*) = \mathbf{0}$ is not a sufficient condition for $\mathbf{J}_h(\mathbf{y}^*) = \mathbf{0}$. Since $S(h(\mathbf{y}_\varepsilon))\geq\mathbf{0}$ as defined in (\ref{eq:prelim_lag}), $h(\mathbf{y}_\varepsilon)=~\mathbf{0}$ can be met if and only if $\mathbf{1}^{\intercal}S(h(\mathbf{y}_\varepsilon)) = 0$ holds. Thus, from Corollary~\ref{cor:positive_error}, $h(\mathbf{y}_\varepsilon) = \mathbf{0}$ is not guaranteed unless $\mathbf{J}_h(\mathbf{y}^*) = \mathbf{0}$. In summary, the feasibility of equality constraints cannot be ensured by solely using Algorithm~\ref{alg:LDF_noW}, but equality embedding is essential as \texttt{DeepLDE} does.

\section{Experiments}
\label{sec:experi}

\subsection{Simulation Setup}
\label{sec:simul_setup}
The proposed \texttt{DeepLDE} is evaluated for three cases as previously conducted by \cite{donti2021DC3}: solving linearly constrained optimization problems with quadratic and non-convex object functions as well as AC-OPF. For a fair comparison with \texttt{DC3}, we test \texttt{DeepLDE} on AC-OPF using the IEEE 57 bus system without branch flow limits (Case 57) and the IEEE 118 bus system with branch flow limits (Case 118). We evaluate the performance in three aspects: \textbf{feasibility} by measuring the maximum and average violation degree of constraints, \textbf{optimality} by calculating the average value of the object function, and the average \textbf{computation time}. All the simulations are conducted on test data.
\begin{table}[t]
	\centering
	\caption{Comparsion with NN-based baselines}
	\label{table:baseline}
	\begin{tabular}{c|ccc}
		\toprule
		\makecell{Method} & \makecell{Self-\\supervised} & \makecell{Equality\\Constraints}  & \makecell{Inequality\\Constraints}\\\midrule\midrule
            \texttt{{SL}}&\xmark&\xmark&\xmark\\
            \texttt{{SL+LD}}\cite{fioretto2020Lagrangian, fioretto2020predicting}&\xmark&\xmark&Primal-dual\\
            \texttt{{SL+E}}\cite{zamzam2020learning,donti2021DC3}&\xmark&Equality embedding&\xmark\\
            \texttt{{SL+LDE}}&\xmark&Equality embedding&Primal-dual\\
            \midrule
            \texttt{{DC3}}\cite{donti2021DC3}&\cmark &Equality embedding&Post-processing\\
            \texttt{{LDF}}\cite{fioretto2020Lagrangian, fioretto2020predicting} &\cmark &\xmark&Primal-dual\\
            \textbf{\texttt{{DeepLDE}}} &\cmark &Equality embedding&Primal-dual\\
		\bottomrule
	\end{tabular}
\end{table}

\begin{figure}[t]
	\centering
	\includegraphics[width=0.99\columnwidth]{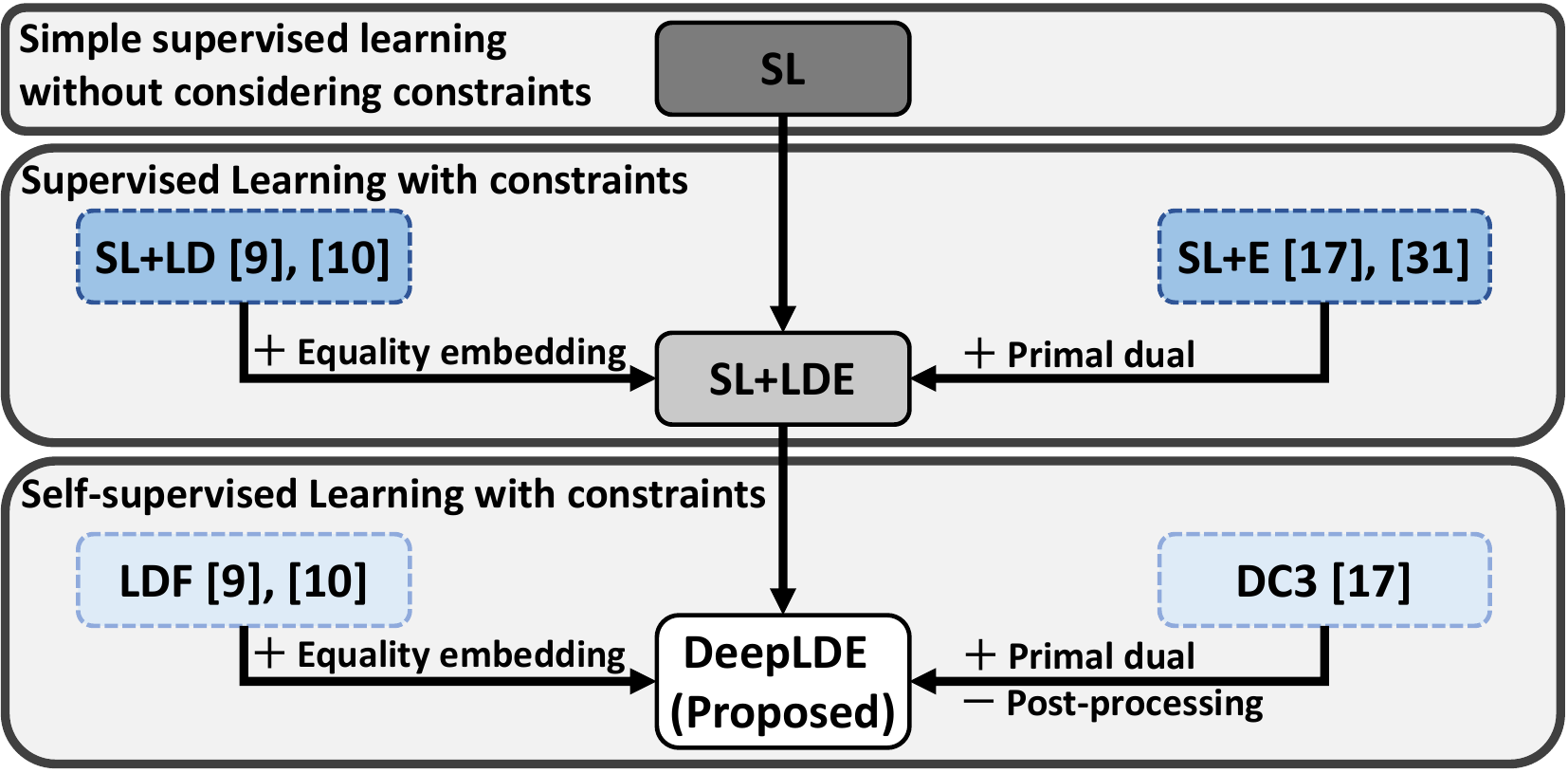}
	\caption{Comparsion of NN-based baselines and \texttt{DeepLDE}.}
	\label{fig:baseline}
\end{figure}
\begin{table}[t]
	\centering
	\caption{The number of iterations of each method.}
	\label{table:iteration}
	\begin{tabular}{c|ccccc}
             \toprule\makecell{Method}&\makecell{$I_w$}&\makecell{$I$}&\makecell{$T$}&\makecell{$\beta$}&\makecell{Total}\\
             \midrule\midrule
            \makecell{\texttt{SL}, \texttt{SL+E}, \texttt{DC3}} & - & - &1,000 & - & 1,000\\
            \midrule
            \makecell{\texttt{SL+LD}, \texttt{SL+LDE}, \texttt{LDF}, \texttt{DeepLDE}} & \makecell{100} & \makecell{25} & \makecell{15} & \makecell{5} & \makecell{1,000}\\
            \bottomrule
	\end{tabular}
\end{table}

We rigorously compare \texttt{{DeepLDE}} with ten state-of-the-art baselines across three categories as follows:

\begin{itemize}
\item \textbf{Conventional solvers without NN}: We use four conventional solvers for solving constrained optimization.
 \begin{description}[leftmargin=!,labelwidth=\widthof{\bfseries \texttt{PYPOWER}}]
 \item[\texttt{OSQP}]: A mathematical solver specialized in solving QP for CPU setting \cite{stellato2020OSQP}.
 \item[\texttt{qpth}]: A differentiable solver specialized in solving QP for GPU setting \cite{amos2017optnet}. 
 \item[\texttt{IPOPT}]: A nonlinear optimization solver based on a primal-dual interior point method \cite{wachter2006implementation}.
 \item[\texttt{PYPOWER}]: An AC-OPF solver which is a port of \texttt{{MATPOWER}} \cite{zimmerman2010matpower} to Python.
 \end{description}
 
\item \textbf{Supervised learning with NN}: We obtain pre-solved solutions from conventional solvers to train NNs in a supervised learning (SL) cases. 
\begin{description}[leftmargin=!,labelwidth=\widthof{\bfseries \texttt{SL+LDE}}]
\item[\texttt{SL}]: Simple SL scheme that minimizes MSE between the predicted output and its corresponding label.
\item[\texttt{SL+LD}]: Lagrange dual framework for \texttt{{SL}} \cite{fioretto2020Lagrangian, fioretto2020predicting}.
\item[\texttt{SL+E}]: Equality embedded SL scheme \cite{zamzam2020learning,donti2021DC3}.
\item[\texttt{SL+LDE}]: SL version of \texttt{DeepLDE}.
\end{description}
  
\item \textbf{Self-supervised learning with NN}: This type of baselines train the NN without relying on any labels. \texttt{DeepLDE} is also in this category.
\begin{description}[leftmargin=!,labelwidth=\widthof{\bfseries \texttt{LDF}}]
\item[\texttt{DC3}]: Use equality embedding and inequality correction to impose constraints to NN \cite{donti2021DC3}.
\item[\texttt{LDF}]: Self-supervised version of  \cite{fioretto2020Lagrangian,fioretto2020predicting} (similar to \texttt{DeepLDE} without equality embedding).
\end{description}
\end{itemize}
Comparisons of NN-based baselines are summarized in Talbe~\ref{table:baseline}. Also, we illustrate how \texttt{DeepLDE} differs from other state-of-the-art methods in Fig.~\ref{fig:baseline}. For a fair comparison, the conventional solvers are assumed to be fully parallelized. 

We train the NN with a \textit{warm-up} period \cite{nandwani2019primal} to enhance the performance of Algorithm~\ref{alg:DeepLDE_noW}. Similarly, all Lagrange dual based approaches, i.e., \texttt{SL+LD}, \texttt{SL+LDF}, and \texttt{LDF}, are also trained with the warm-up. In the warm-up period, the NN is trained with an additional inner iteration before the first outer iteration. This additional process does not harm the convergence of the learning process and reduces the optimality gap if it exists. The maximum warm-up iteration $I_w$, inner iteration $I$, and outer iteration $T$, as well as increments of outer iteration $\beta$ are described in Table~\ref{table:iteration}. Additionally, the total number of iterations is determined by $I_w + \sum_{t=1}^T I + \beta (t-1)$. For a fair comparison, we obtain the same total number of iterations across all methods using $T$, $I_w$, $I$, and $\beta$ as shown in Table~\ref{table:iteration}. 
\begin{table}[t]
	\centering
	\caption{Hyperparameters of simulations.}
	\label{table:QP_hyper}
	\begin{tabular}{c|c|cccc}
		\toprule
		\multicolumn{2}{c|}{\makecell{Method}} & $\rho_0$ & $s_0$  & $\lambda_0$  & $\mu_0$\\\midrule\midrule
            \multirow{4}{*}{QP}&\texttt{\makecell{SL+LD,\\LDF}}&\makecell{\textbf{0.1},\\0.5}&\makecell{0.1,\\\textbf{0.5}}&\makecell{0.0,\\\textbf{0.1}}&\makecell{0.0,\\\textbf{0.1}}\\\cmidrule(lr){2-6}
            &\texttt{\makecell{SL+LDE,\\DeepLDE}}&\makecell{\textbf{0.1},\\0.5}&-&\makecell{0.0,\\\textbf{0.1}}&-\\\midrule
            \multirow{4}{*}{\makecell{Non-\\convex}}&\texttt{\makecell{SL+LD,\\ LDF}}&\makecell{\textbf{0.0001},\\0.0005}&\makecell{0.0001,\\\textbf{0.0005}}&\makecell{0.0,\\\textbf{0.1}}&\makecell{0.0,\\\textbf{0.1}}\\\cmidrule(lr){2-6}
            &\texttt{\makecell{SL+LDE,\\DeepLDE}}&\makecell{\textbf{0.0001},\\0.0005}&-&\makecell{0.0,\\\textbf{0.1}}&-\\\midrule
            \multirow{4}{*}{Case 57}&\texttt{\makecell{SL+LD,\\LDF}}&\makecell{\textbf{0.1},\\0.5}&\makecell{0.1,\\\textbf{0.5}}&\makecell{0.0,\\\textbf{0.1}}&\makecell{0.0,\\\textbf{0.1}}\\\cmidrule(lr){2-6}
            &\texttt{\makecell{SL+LDE,\\DeepLDE}}&\makecell{\textbf{0.1},\\0.5}&-&\makecell{0.0,\\\textbf{0.1}}&-\\\midrule
            \multirow{4}{*}{Case 118}&\texttt{\makecell{SL+LD,\\LDF}}&\makecell{\textbf{0.01},\\0.05}&\makecell{0.01,\\\textbf{0.05}}&\makecell{0.0,\\\textbf{0.1}}&\makecell{0.0,\\\textbf{0.1}}\\\cmidrule(lr){2-6}
            &\texttt{\makecell{SL+LDE,\\DeepLDE}}&\makecell{\textbf{0.01},\\0.05}&-&\makecell{0.0,\\\textbf{0.1}}&-\\
		\bottomrule
	\end{tabular}
\end{table}
In Tables~\ref{table:QP_hyper}, we list the candidates of hyperparameters of entire simulations. The selected values are in \textbf{bold}. Also, we set the learning rate $\eta$ to $10^{-3}$. For \texttt{DC3}, we use the same hyperparameters selected by \cite{donti2021DC3}. However, we found that using the reported step size of correction above $10^{-4}$ for \texttt{DC3} in Case 118 causes instability during training. Hence, we use $10^{-4}$ of step size of correction for \texttt{DC3} in Case 118. The rest of the hyperparameters of \texttt{DC3} in Case 118 are all equal to Case 57.

For all methods and simulations, we set the minibatch size to 200 and $\gamma$ to 0.01. We use the NN with two hidden layers with 200 neurons, including the ELU activation function, and a dropout rate of 0.1 for all experiments \cite{clevert2015fast, srivastava2014dropout}. We use Adam optimizer \cite{kingma2014adam} to train NN. To obtain statistically reliable results, we repeat experiments five times as done in \cite{donti2021DC3} and measure the standard deviations, which are provided in the parenthesis right below the average value, see Table~\ref{table:QP}$-$\ref{table:ACOPF}. We find that, in most cases, the standard deviations are nearly zero, which confirms the validity of our experiments. All the experiments are conducted on NVIDIA A100 SXM4 40GB GPU and Intel Xeon 2.20GHz CPU.

\begin{table*}[ht]
	\centering
	\caption{Results on linearly constrained QP.}
	\label{table:QP}
	\begin{tabular}{c|cccccc||ccccccc}
 \toprule
            \multirow{3}{*}{\makecell{Method}}&\multicolumn{6}{c||}{$n_\text{eq} = 70, n_\text{ineq} = 30$}&\multicolumn{6}{c}{$n_\text{eq} = 30, n_\text{ineq} = 70$}\\
            \cmidrule(lr){2-13}
		&\makecell{Max\\Eq.}  & \makecell{Mean\\Eq.}  & \makecell{Max\\Ineq.} & \makecell{Mean\\Ineq.} & \makecell{Obj.\\value} &\makecell{Time\\(sec)} & \makecell{Max\\Eq.}  & \makecell{Mean\\Eq.}  & \makecell{Max\\Ineq.} & \makecell{Mean\\Ineq.} & \makecell{Obj.\\value}&\makecell{Time\\(sec)}\\\midrule\midrule
            \texttt{{OSQP}}  & \makecell{0.00\\(0.00)}& \makecell{0.00\\(0.00)} &\makecell{0.00\\(0.00)}&\makecell{0.00\\(0.00)}& \makecell{-14.88\\(0.00)} & \makecell{0.002\\(0.00)} & \makecell{0.00\\(0.00)}& \makecell{0.00\\(0.00)} &\makecell{0.00\\(0.00)}&\makecell{0.00\\(0.00)} & \makecell{-21.01\\(0.00)}& \makecell{0.002\\(0.00)} \\
            \texttt{{qpth}} & \makecell{0.00\\(0.00)}& \makecell{0.00\\(0.00)} &\makecell{0.00\\(0.00)}&\makecell{0.00\\(0.00)} & \makecell{-14.88\\(0.00)}& \makecell{2.261\\(0.344)} & \makecell{0.00\\(0.00)}& \makecell{0.00\\(0.00)} &\makecell{0.00\\(0.00)}&\makecell{0.00\\(0.00)} & \makecell{-21.01\\(0.00)}& \makecell{2.108\\(0.017)}\\
            \midrule
            \texttt{{SL}}&\textcolor{red}{\makecell{{0.21}\\(0.01)}}&\textcolor{red}{\makecell{{0.02}\\(0.00)}}&\makecell{{0.03}\\(0.00)}&\makecell{{0.00}\\(0.00)} &\makecell{-14.90\\(0.01)} &\makecell{0.005\\0.000}&\textcolor{red}{\makecell{{0.03}\\(0.01)}}&{\makecell{{0.00}\\(0.00)}}&\textcolor{red}{\makecell{{0.09}\\(0.00)}}&\textcolor{red}{\makecell{{0.03}\\(0.00)}} &\makecell{-21.02\\(0.00)}&\makecell{0.005\\0.000}\\
            \texttt{{SL+LD}}&\textcolor{red}{\makecell{{0.08}\\(0.01)}}&\textcolor{red}{\makecell{{0.01}\\(0.00)}}&\makecell{{0.00}\\(0.00)}&\makecell{{0.00}\\(0.00)}&\makecell{-14.79\\(0.04)}&{\makecell{{0.005}\\(0.000)}}&\textcolor{red}{\makecell{{0.09}\\(0.00)}}&\textcolor{red}{\makecell{{0.03}\\(0.00)}}&\makecell{{0.00}\\(0.00)}&\makecell{{0.00}\\(0.00)}&\makecell{-20.76\\(0.00)}&{\makecell{{0.005}\\(0.000)}}\\
            \texttt{{SL+E}}&\makecell{{0.00}\\(0.00)} & \makecell{{0.00}\\(0.00)}&\textcolor{red}{\textcolor{red}{\makecell{{0.23}\\(0.01)}}}&\textcolor{red}{\textcolor{red}{\makecell{{0.01}\\(0.00)}}}&\makecell{-14.87\\(0.01)}&\makecell{0.008\\(0.000)}&\makecell{{0.00}\\(0.00)} & \makecell{{0.00}\\(0.00)}&\textcolor{red}{\textcolor{red}{\makecell{{0.31}\\(0.03)}}}&\textcolor{red}{\textcolor{red}{\makecell{{0.02}\\(0.00)}}}&\makecell{-21.00\\(0.02)}&\makecell{0.008\\(0.000)}\\
            \texttt{{SL+LDE}}& \makecell{{0.00}\\(0.00)} & \makecell{{0.00}\\(0.00)} & \makecell{{0.00}\\(0.00)} &\makecell{{0.00}\\(0.00)}&\makecell{-14.60\\(0.02)}&\makecell{0.008\\(0.000)}& \makecell{{0.00}\\(0.00)} & \makecell{{0.00}\\(0.00)} & \makecell{{0.00}\\(0.00)} &\makecell{{0.00}\\(0.00)}&\makecell{-19.87\\(0.02)}&\makecell{0.008\\(0.000)}\\
            \midrule
            \multicolumn{13}{c}{\quad\quad\quad\quad\quad\quad\quad Self-supervised learning (without labels)}\\
            \midrule
            \texttt{{DC3}} & \makecell{{0.00}\\(0.00)} & \makecell{{0.00}\\(0.00)} & \makecell{{0.00}\\(0.00)} &\makecell{{0.00}\\(0.00)} &\makecell{-14.47\\(0.01)}& \makecell{0.008\\(0.000)}& \makecell{0.00\\(0.00)} & \makecell{{0.00}\\(0.00)} & \makecell{{0.00}\\(0.00)} &\makecell{{0.00}\\(0.00)} &\makecell{-19.66\\(0.01)}& \makecell{0.008\\(0.000)}\\
            \texttt{{LDF}} & \textcolor{red}{\makecell{{0.22}\\(0.00)}} & \textcolor{red}{\makecell{{0.08}\\(0.00)}} & \makecell{{0.00}\\(0.00)} &\makecell{{0.00}\\(0.00)} &\makecell{-14.73\\(0.00)} & \makecell{0.005\\(0.000)}& \textcolor{red}{\makecell{{0.19}\\(0.00)}} & \textcolor{red}{\makecell{{0.06}\\(0.00)}} & \makecell{{0.00}\\(0.00)} &\makecell{{0.00}\\(0.00)} &\makecell{-20.91\\(0.00)} & \makecell{0.005\\(0.000)}\\
            \textbf{\texttt{{DeepLDE}}}&\textbf{\makecell{{0.00}\\(0.00)}} & \textbf{\makecell{{0.00}\\(0.00)}} & \textbf{\makecell{{0.00}\\(0.00)}}& \textbf{\makecell{{0.00}\\(0.00)}}  & \makecell{\textbf{-14.67}\\\textbf{(0.00)}}& \textbf{\makecell{0.008\\(0.000)}}&\textbf{\makecell{{0.00}\\(0.00)}} & \textbf{\makecell{{0.00}\\(0.00)}} & \textbf{\makecell{{0.00}\\(0.19)}}& \textbf{\makecell{{0.00}\\(0.00)}} & \makecell{\textbf{-19.97}\\\textbf{(0.40)}}& \textbf{\makecell{0.008\\(0.000)}}\\
		\bottomrule
	\end{tabular}
\end{table*}

\begin{table*}[ht]
	\centering
	\caption{Results on linearly constrained non-convex programming.}
	\label{table:QP-v}
	\begin{tabular}{c|cccccc||ccccccc}
 \toprule
            \multirow{4}{*}{\makecell{Method}}&\multicolumn{6}{c||}{$n_\text{eq} = 70, n_\text{ineq} = 30$}&\multicolumn{6}{c}{$n_\text{eq} = 30, n_\text{ineq} = 70$}\\
            \cmidrule(lr){2-13}  & \makecell{Max\\Eq.}  & \makecell{Mean\\Eq.}  & \makecell{Max\\Ineq.} & \makecell{Mean\\Ineq.}& \makecell{Obj.\\value} &\makecell{Time\\(sec)} & \makecell{Max\\Eq.}  & \makecell{Mean\\Eq.}  & \makecell{Max\\Ineq.} & \makecell{Mean\\Ineq.} & \makecell{Obj.\\value}&\makecell{Time\\(sec)}\\\midrule\midrule
            \makecell{\texttt{{IPOPT}}}  &  \makecell{0.00\\(0.00)} & \makecell{0.00\\0.00} &\makecell{0.00\\(0.00)}&\makecell{0.00\\(0.00)} & \makecell{-9.93\\(0.00)}& \makecell{0.207\\(0.012)}& \makecell{0.00\\(0.00)} & \makecell{0.00\\(0.00)} &\makecell{0.00\\(0.00)}&\makecell{0.00\\(0.00)} & \makecell{-15.65\\(0.00)}& \makecell{0.270\\(0.000)} \\
            \midrule
            \texttt{{SL}}&\textcolor{red}{\makecell{{0.24}\\(0.01)}}&\textcolor{red}{\makecell{{0.08}\\(0.00)}}&\textcolor{red}{\makecell{{0.03}\\(0.00)}}&\makecell{{0.00}\\(0.00)} &\makecell{-9.96\\(0.00)}&\makecell{0.005\\0.000}&\textcolor{red}{\makecell{{0.09}\\(0.01)}}&\textcolor{red}{\makecell{{0.03}\\(0.00)}}&\textcolor{red}{\makecell{{0.04}\\(0.00)}}&\makecell{{0.00}\\(0.00)} &\makecell{-15.66\\(0.00)}&\makecell{0.005\\(0.000)}\\
            \texttt{{SL+LD}}&\textcolor{red}{\makecell{{0.23}\\(0.02)}}&\textcolor{red}{\makecell{{0.08}\\(0.01)}}&\makecell{{0.00}\\(0.00)}&\makecell{{0.00}\\(0.00)}&\makecell{-9.95\\(0.01)}&\textcolor{red}{\makecell{{0.005}\\(0.000)}}&\textcolor{red}{\makecell{{0.09}\\(0.00)}}&\textcolor{red}{\makecell{{0.03}\\(0.00)}}&\makecell{{0.00}\\(0.00)}&\makecell{{0.00}\\(0.00)}&\makecell{-15.62\\(0.00)}&{\makecell{{0.005}\\(0.000)}}\\
            \texttt{{SL+E}}&\makecell{{0.00}\\(0.00)} & \makecell{{0.00}\\(0.00)}&\textcolor{red}{\textcolor{red}{\makecell{{0.13}\\(0.04)}}}&\textcolor{red}{\textcolor{red}{\makecell{{0.01}\\(0.00)}}}&\makecell{-9.90\\(0.01)}&\makecell{0.008\\(0.000)}&\makecell{{0.00}\\(0.00)} & \makecell{{0.00}\\(0.00)}&\textcolor{red}{\textcolor{red}{\makecell{{0.23}\\(0.01)}}}&\textcolor{red}{\textcolor{red}{\makecell{{0.02}\\(0.00)}}}&\makecell{-15.64\\(0.00)}&\makecell{0.008\\(0.000)}\\
            \texttt{{SL+LDE}}& \makecell{{0.00}\\(0.00)} & \makecell{{0.00}\\(0.00)} & \makecell{{0.00}\\(0.00)} &\makecell{{0.00}\\(0.00)}&\makecell{-9.73\\(0.05)}&\makecell{0.008\\(0.000)}& \makecell{{0.00}\\(0.00)} & \makecell{{0.00}\\(0.00)} & \makecell{{0.00}\\(0.00)} &\makecell{{0.00}\\(0.00)}&\makecell{-14.96\\(0.01)}&\makecell{0.008\\(0.000)}\\
            \midrule
            \multicolumn{13}{c}{\quad\quad\quad\quad\quad\quad\quad Self-supervised learning (without labels)}\\
            \midrule
            \texttt{{DC3}} & \makecell{{0.00}\\(0.00)} & \makecell{{0.00}\\(0.00)} & \makecell{{0.00}\\(0.00)} &\makecell{{0.00}\\(0.00)} &\makecell{-9.84\\(0.01)}& \makecell{0.008\\(0.000)}& \makecell{{0.00}\\(0.00)} & \makecell{{0.00}\\(0.00)} & \makecell{{0.00}\\(0.00)} &\makecell{{0.00}\\(0.00)} &\makecell{-14.76\\(0.02)}& \makecell{0.008\\(0.005)}\\
            \texttt{{LDF}} & \textcolor{red}{\makecell{{0.22}\\(0.00)}} & \textcolor{red}{\makecell{{0.08}\\(0.00)}} & \makecell{{0.00}\\(0.00)} &\makecell{{0.00}\\(0.00)} &\makecell{-9.96\\(0.00)}& \makecell{0.005\\(0.000)}& \textcolor{red}{\makecell{{0.11}\\(0.03)}} & \textcolor{red}{\makecell{{0.03}\\(0.00)}} & \makecell{{0.00}\\(0.00)} &\makecell{{0.00}\\(0.00)} &\makecell{-15.56\\(0.00)}& \makecell{0.005\\(0.000)}\\
            \textbf{\texttt{{DeepLDE}}} &\textbf{\makecell{{0.00}\\(0.00)}} & \textbf{\makecell{{0.00}\\(0.00)}} & \textbf{\makecell{{0.00}\\(0.00)}}& \textbf{\makecell{{0.00}\\(0.00)}} & \makecell{\textbf{-9.89}\\\textbf{(0.00)}}& \textbf{\makecell{0.008\\(0.000)}}&\textbf{\makecell{{0.00}\\(0.00)}} & \textbf{\makecell{{0.00}\\(0.00)}} & \textbf{\makecell{{0.00}\\(0.00)}}& \textbf{\makecell{{0.00}\\(0.00)}} & \makecell{\textbf{-15.15}\\\textbf{(0.02)}}& \textbf{\makecell{0.008\\(0.000)}}\\
		\bottomrule
	\end{tabular}
\end{table*}
\subsection{Linearly Constrained Nonlinear Programming}
\label{sec:exp-linear}
For performance evaluation, we first consider the following linearly constrained quadratic programming \cite{donti2021DC3}:
\begin{equation}
\begin{aligned}
    &\min_{\mathbf{y}\in\mathbb{R}^n} \dfrac{1}{2}\mathbf{y}^{\intercal}\mathbf{Q}\mathbf{y} + \mathbf{p}^{\intercal}\mathbf{y},\\
    &\textcolor{black}{\text{s.t.}} \quad \mathbf{A}\mathbf{y} = \mathbf{d},\quad \mathbf{G}\mathbf{y} \leq \mathbf{h},
    \label{eq:QP}
\end{aligned}
\end{equation}
where $\mathbf{Q}\in S_{+}^{n}, \mathbf{p}\in \mathbb{R}^n, \mathbf{A}\in\mathbb{R}^{n_{\text{eq}}\times n}, \mathbf{G}\in\mathbb{R}^{n_{\text{ineq}}
\times n}$ and $\mathbf{h}\in\mathbb{R}^{n_{\text{ineq}}}$. The input variable vector $\mathbf{d}\in\mathbb{R}^{n_{\text{eq}}}$ is sampled from $\mathcal{D}$.
The data generation process is provided in Appendix~\ref{appendix:data-generate}. We also consider a non-convex objective function with the same constraints of (\ref{eq:QP}) as follows:
\begin{equation}
\begin{aligned}
    &\min_{\mathbf{y}\in\mathbb{R}^n} \dfrac{1}{2}\mathbf{y}^{\intercal}\mathbf{Q}\mathbf{y} + \mathbf{p}^{\intercal}\textbf{sin}(\mathbf{y}),\\
    &\textcolor{black}{\text{s.t.}} \quad \mathbf{A}\mathbf{y} = \mathbf{d},\quad \mathbf{G}\mathbf{y} \leq \mathbf{h},
    \label{eq:NonConv}
\end{aligned}
\end{equation}
 where $\textbf{sin}(\cdot)$ is the element-wise sine function. We set $n_{\text{eq}} = 70$, $n_{\text{ineq}} = 30$ and $n_{\text{eq}} = 30$, $n_{\text{ineq}} = 70$ for comparison \cite{donti2021DC3}.
 \subsubsection{Feasibility}
Tables \ref{table:QP} and \ref{table:QP-v} present the results of conventional optimization solvers and NN-based approaches in solving (\ref{eq:QP}) and (\ref{eq:NonConv}), respectively. We use \textcolor{red}{red} letters to highlight the degree of violations of NN-based methods. Recall that equality embedding is for ensuring equality constraints and the primal-dual method is for ensuring inequality constraints. As can be seen, the equality embedding method used in \texttt{{SL+LDE}}, \texttt{{DC3}}, and \texttt{{DeepLDE}} always satisfies the equality constraints. Similarly, the primal-dual method used in \texttt{SL+LD}, \texttt{SL+LDE}, \texttt{LDF}, and \texttt{DeepLDE} effectively enforces the inequality constraints. By contrast, \texttt{SL} is the poorest in satisfying the equality and the inequality constraints because it uses neither equality embedding nor the primal-dual method. Similarly, \texttt{SL+E} violates the inequality constraints, and \texttt{SL+LD} and \texttt{LDF} violate the equality constraints. By contrast, \texttt{{SL+LDE}}, \texttt{{DC3}}, and \texttt{{DeepLDE}} yield feasible solutions because of using the equality embedding and/or the primal-dual method. Then, the next question is how good (optimality) and how fast (computation time) the algorithms are.
\subsubsection{Optimality}
While \texttt{SL+E} achieves the smallest optimality gap among the NN-based methods, its solutions are infeasible, and thus meaningless. Among the NN methods that obtain \textit{feasible} solutions (\texttt{SL+LDE}, \texttt{DC3}, and \texttt{DeepLDE}), \texttt{DeepLDE} achieves the smallest optimality gap (e.g., 1.41\% in QP and 0.60\% in non-convex programming).
\subsubsection{Computation time}
As shown in Tables~\ref{table:QP} and \ref{table:QP-v}, \texttt{DC3} exhibits the longest computation time among the NN-based approaches due to the additional iterations required in the correction process. It is worth noting that, although \texttt{DeepLDE} does not use pre-solved labels for training, its computation time is similar to \texttt{SL+LDE} while achieving smaller optimality gap. In solving QP, the proposed \texttt{DeepLDE} is 458 times faster than \texttt{qpth}; \texttt{OSQP} is the fastest but cannot solve other non-convex problems. By contrast, \texttt{DeepLDE} is a general-purpose solver and 28 times faster than the off-the-shelf general-purpose solver \texttt{IPOPT} for solving non-convex optimization, as shown in Table \ref{table:QP-v}.

\begin{table*}[ht]
	\centering
	\caption{Results on AC-OPF.}
	\label{table:ACOPF}
	\begin{tabular}{c|cccccc||ccccccc}
 \toprule
            \multirow{3}{*}{\makecell{Method}}&\multicolumn{6}{c||}{Case 57 (without branch flow limits)}&\multicolumn{6}{c}{Case 118 (with branch flow limits)}\\
            \cmidrule(lr){2-13}
            & \makecell{Max\\Eq.}  & \makecell{Mean\\Eq.}  & \makecell{Max\\Ineq.} & \makecell{Mean\\Ineq.} & \makecell{Obj.\\value} &\makecell{Time\\(sec)} & \makecell{Max\\Eq.}  & \makecell{Mean\\Eq.}  & \makecell{Max\\Ineq.} & \makecell{Mean\\Ineq.} & \makecell{Obj.\\value} &\makecell{Time\\(sec)}\\\midrule\midrule
            \makecell{\texttt{{PYPOWER}}}   & \makecell{{0.00}\\(0.00)} & \makecell{{0.00}\\(0.00)} &\makecell{{0.00}\\(0.00)}&\makecell{{0.00}\\(0.00)} & \makecell{3.74\\(0.00)} & \makecell{0.690\\(0.000)} & \makecell{{0.00}\\(0.00)} & \makecell{{0.00}\\(0.00)} &\makecell{{0.00}\\(0.00)}&\makecell{{0.00}\\(0.00)} & \makecell{12.91\\(0.00)}& \makecell{1.063\\(0.000)}\\
            \midrule
            \texttt{{SL}}&\textcolor{red}{\makecell{{0.21}\\(0.01)}}&\textcolor{red}{\makecell{{0.02}\\(0.00)}}&\makecell{{0.00}\\(0.00)}&\makecell{{0.00}\\(0.00)} &\makecell{3.72\\(0.02)}&\makecell{0.017\\0.000}&\textcolor{red}{\makecell{{0.45}\\(0.02)}}&\textcolor{red}{\makecell{{0.06}\\(0.00)}}&\makecell{{0.00}\\(0.00)}&\makecell{{0.00}\\(0.00)} &\makecell{12.94\\(0.06)}&\makecell{0.017\\(0.000)}\\
            \texttt{{SL+LD}}&\textcolor{red}{\makecell{{0.10}\\(0.01)}}&\textcolor{red}{\makecell{{0.01}\\(0.00)}}&\makecell{{0.00}\\(0.00)}&\makecell{{0.00}\\(0.00)}&\makecell{3.75\\(0.02)}&{\makecell{{0.018}\\(0.001)}}&\textcolor{red}{\makecell{{0.36}\\(0.02)}}&\textcolor{red}{\makecell{{0.04}\\(0.00)}}&\makecell{{0.00}\\(0.00)}&\makecell{{0.00}\\(0.00)}&\makecell{13.14\\(0.03)}&{\makecell{{0.017}\\(0.001)}}\\
            \texttt{{SL+E}}&\makecell{{0.00}\\(0.00)} & \makecell{{0.00}\\(0.00)}&\textcolor{red}{\textcolor{red}{\makecell{{0.04}\\(0.00)}}}&{{\makecell{{0.00}\\(0.00)}}}&\makecell{3.74\\(0.00)}&\makecell{0.039\\(0.001)}&\makecell{{0.00}\\(0.00)} & \makecell{0.00\\(0.00)}&\textcolor{red}{\textcolor{red}{\makecell{0.10\\(0.00)}}}&{\textcolor{red}{\makecell{{0.00}\\(0.00)}}}&\makecell{13.25\\(0.00)}&\makecell{0.053\\(0.002)}\\
            \texttt{{SL+LDE}}& \makecell{{0.00}\\(0.00)} & \makecell{{0.00}\\(0.00)} & \makecell{{0.00}\\(0.00)} &\makecell{{0.00}\\(0.00)}&\makecell{3.74\\(0.00)}&\makecell{0.039\\(0.001)}& \makecell{{0.00}\\(0.00)} & \makecell{{0.00}\\(0.00)} & \makecell{{0.00}\\(0.00)} &\makecell{{0.00}\\(0.00)}&\makecell{13.27\\(0.01)}&\makecell{0.051\\(0.000)}\\
            \midrule
            \multicolumn{13}{c}{\quad\quad\quad\quad\quad\quad\quad Self-supervised learning (without labels)}\\
            \midrule
            \makecell{\texttt{{DC3}}} & \makecell{{0.00}\\(0.00)} & \makecell{{0.00}\\(0.00)} & \makecell{{0.00}\\(0.00)} &\makecell{{0.00}\\(0.00)} &\makecell{3.74\\(0.00)}& \makecell{0.066\\(0.000)}& \makecell{{0.00}\\(0.00)} & \makecell{{0.00}\\(0.00)} & \textcolor{red}{\makecell{{0.01}\\(0.00)}} &\makecell{{0.00}\\(0.00)} &\makecell{13.32\\(0.08)}& \makecell{0.284\\(0.002)}\\
            \makecell{\texttt{{LDF}}} & \textcolor{red}{\makecell{{0.10}\\(0.00)}} & \textcolor{red}{\makecell{{0.02}\\(0.00)}} & \makecell{{0.00}\\(0.00)} &\makecell{{0.00}\\(0.00)} &\makecell{3.76\\(0.00)}& \makecell{0.017\\(0.000)}& \textcolor{red}{\makecell{{0.33}\\(0.03)}} & \textcolor{red}{\makecell{{0.04}\\(0.00)}} & \makecell{{0.00}\\(0.00)} &\makecell{{0.00}\\(0.00)} &\makecell{14.74\\(0.12)}& \makecell{0.017\\(0.000)}\\
            \textbf{\texttt{{DeepLDE}}} &\textbf{\makecell{{0.00}\\(0.00)}} & \textbf{\makecell{{0.00}\\(0.00)}} & \textbf{\makecell{{0.00}\\(0.00)}}& \textbf{\makecell{{0.00}\\(0.00)}} & \makecell{\textbf{3.74}\\\textbf{(0.00)}}& \textbf{\makecell{0.039\\(0.001)}}&\textbf{\makecell{{0.00}\\(0.00)}} & \textbf{\makecell{{0.00}\\(0.00)}} & \textbf{\makecell{{0.00}\\(0.00)}}& \textbf{\makecell{{0.00}\\(0.00)}} & \makecell{\textbf{13.18}\\\textbf{(0.00)}}& \textbf{\makecell{0.050\\(0.002)}}\\
		\bottomrule
	\end{tabular}
\end{table*}

\subsection{AC Optimal Power Flow}
\label{sec:AC-OPF}
\begin{figure}[t]
	\centering
	\includegraphics[width=0.95\columnwidth]{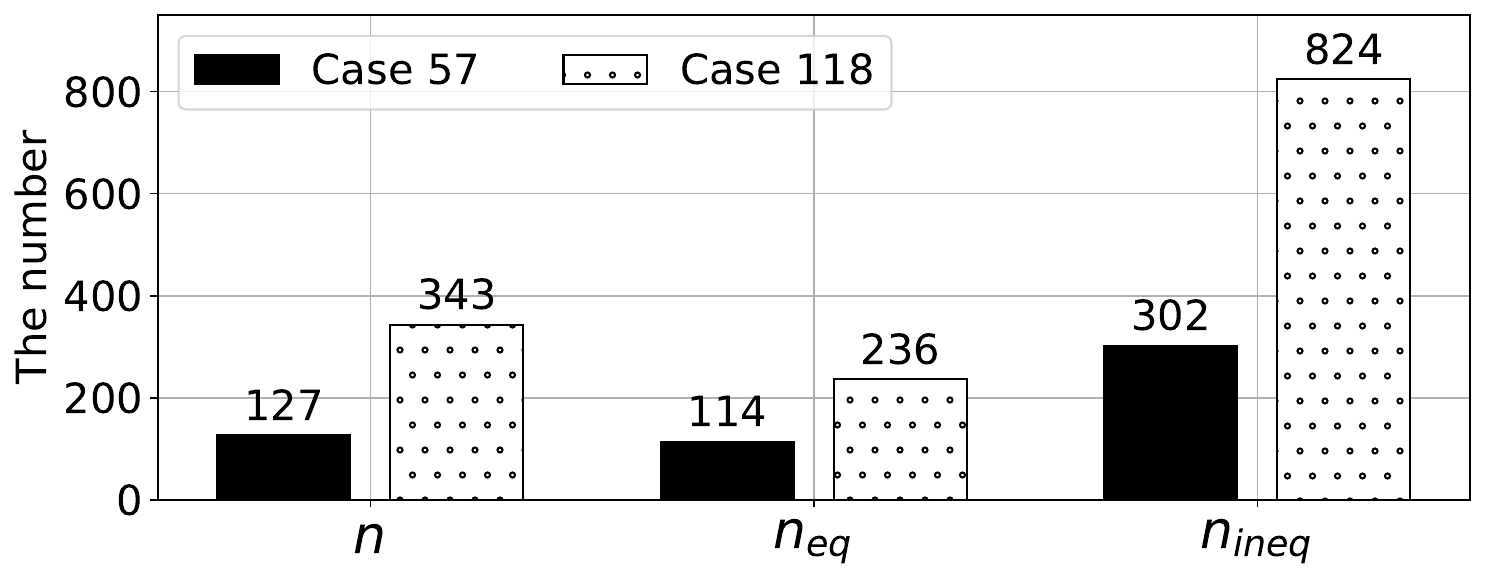}
	\caption{The number of optimization variables ($n$) and constraints ($n_{\text{eq}}$ and $n_{\text{ineq}}$) in Case 57 (black) and Case 118 (dotted).}
	\label{fig:number_opt}
\end{figure}
AC-OPF, one of the most fundamental problems in power systems, is a large-scale real-world optimization problem that minimizes overall power generation cost under non-convex physical constraints. AC-OPF considers a power network (graph) consisting of buses (or nodes) $\mathcal{N}$ and lines (or edges)  $\mathcal{E}$ as follows \cite{pan2022deepopf}:

\begin{subequations}\label{eq:AC-OPF}
    \begin{align}
&\min_{\substack{p_i^g, q_i^g,\\v_i, \forall i\in \mathcal{N}}}&&\sum_{i\in \mathcal{N}} C_i(p_i^g)\label{eq:generation_cost}&\\ 
&\quad\text{s.t.}&&\underline{p}_i^g\leq p_i^g\leq\overline{p}_i^g, \forall i\in\mathcal{N},\label{eq:generator_active_limit}&\\
&\text{}&&\underline{q}_i^g \leq q_i^g\leq\overline{q}_i^g,\label{eq:generator_reactive_limit}&\\
&\text{}
&&|\underline{v}_i|\leq |v_i|\leq|\overline{v}_i|,\label{eq:voltage_magnitude_limit}&\\
&\text{}&&|v_i(v_i^* - v_j^*)Y^*_{ij}|\leq \overline{S}_{ij}, \forall(i,j)\in \mathcal{E}\label{eq:branch_flow_limit},&\\
&\text{}&&\begin{aligned}
    p_i^g-p_i^d =\sum_{(i,j)\in\mathcal{E}}\text{Re}(v_i(v_i^*-v_j^*)Y_{ij}^*)\label{eq:active_pf},
\end{aligned}\\
&\text{}&&\begin{aligned}
    q_i^g-q_i^d =\sum_{(i,j)\in\mathcal{E}}\text{Im}(v_i(v_i^*-v_j^*)Y_{ij}^*)\label{eq:reactive_pf},
\end{aligned}
\end{align}
\end{subequations}

where $C_i(\cdot)$ is the generation cost function. $p_i^g\in\mathbb{R}$ and $q_i^g\in\mathbb{R}$ are active and reactive power generation of bus $i$ and restricted by (\ref{eq:generator_active_limit})$-$(\ref{eq:generator_reactive_limit}). Likewise, $|v_i|$ is bounded by (\ref{eq:voltage_magnitude_limit}) where $v_i\in\mathbb{C}$ is the complex voltage of bus $i$. The maximum transmission power is constrained by the branch flow limits (\ref{eq:branch_flow_limit}). Note that overline and underline in (\ref{eq:generator_active_limit})$-$(\ref{eq:branch_flow_limit}) denote the upper and lower bounds of the inequalities. Finally, $p_i^g$ and $q_i^g$ are derived by power flow equations (\ref{eq:active_pf})$-$(\ref{eq:reactive_pf}) where $p_i^d\in\mathbb{R}$ and $q_i^d\in\mathbb{R}$ are active and reactive power demand of bus $i$, and $Y_{ij}\in\mathbb{C}$ denote the admittance of line $(i,j)$. 

The number of optimization variables is described in Fig.~\ref{fig:number_opt}, which reflects the complexity of two cases quantitatively. All the NN-based methods take $p_i^d$ and $q_i^d$ as an input (i.e., $\mathbf{d}$ in (\ref{eq:prelim_opt})). The output of the NN-based methods with equality embedding is $p_i^g$ and $|v_i|$ (i.e., $\mathbf{x}$ in (\ref{eq:method_reformulZ})), and the rest of the variables are determined by solving (\ref{eq:active_pf})$-$(\ref{eq:reactive_pf}) using $\Psi_h$ as shown in (\ref{eq:method_separateY}). The data generation process is provided in Appendix~\ref{appendix:data-generate}. Table~\ref{table:ACOPF} presents the results in three aspects: feasibility, optimality, and computation time.
\subsubsection{Feasibility}
Since \texttt{DeepLDE} uses both the equality embedding and the primal-dual method, no constraints are violated in both Case 57 and Case 118. By contrast, although \texttt{DC3} uses the equality embedding method, the equality constraints turn out to be violated in Case 118 after using the correction method to satisfy the inequality constraints, see Table~\ref{table:ACOPF}. This infeasible result can be critical because in the case of a large-scale power system optimization problem such as Case 118, attempting to satisfy one kind of constraints (e.g., inequality) may violate the readily satisfied constraints (e.g., equality). This result implies that the correction method used in \texttt{DC3} is less effective in ensuring feasibility for large-scale optimization problems than the primal-dual method of \texttt{DeepLDE}.
\begin{figure}[t]
	\centering
	\includegraphics[width=0.99\columnwidth]{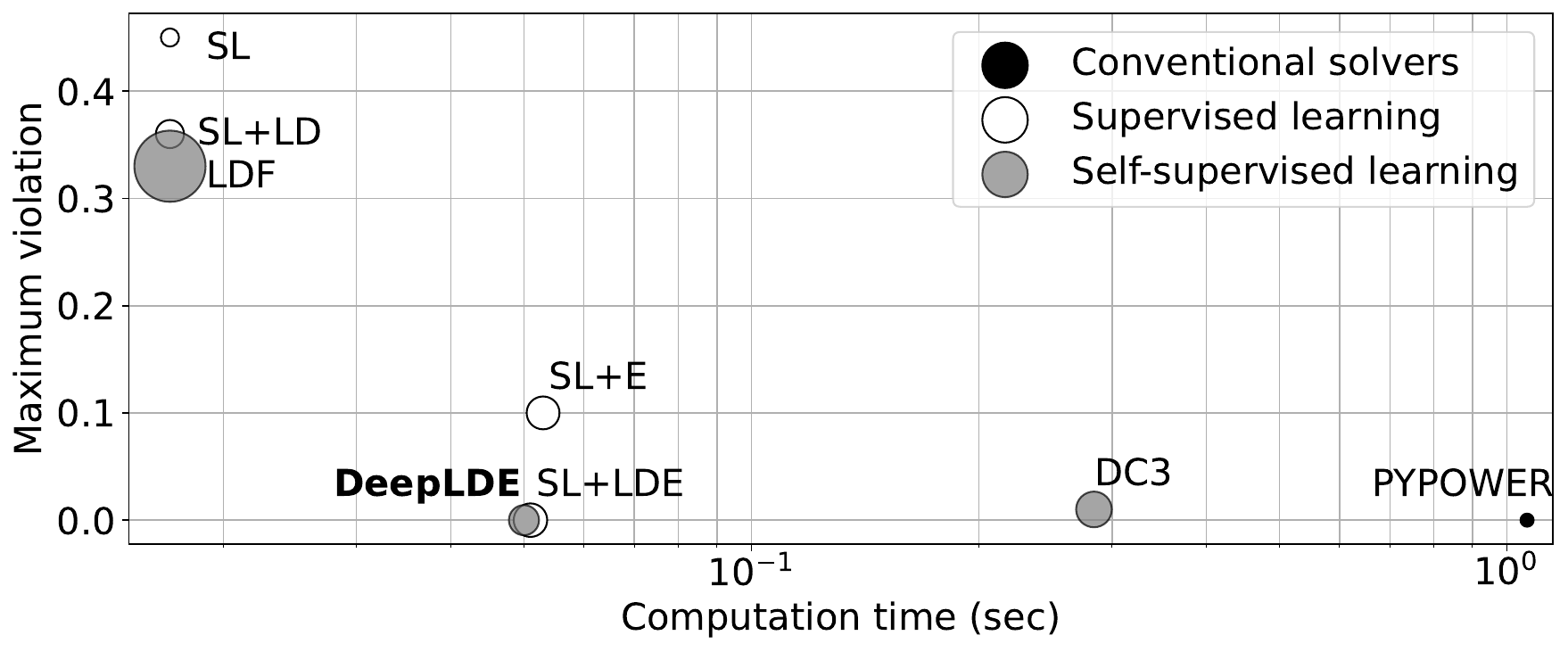}
	\caption{Tradeoff between computation time (x-axis, log scale) and maximum violation (y-axis) in Case 118. The size of the circle is proportional to the optimality gap.}
	\label{fig:case118_tradeoff}
\end{figure}
\subsubsection{Optimality}
Among all methods, the largest optimality gap is observed in \texttt{LDF} which solely uses the primal-dual method without using labels. This gap further increases in the number of buses. Furthermore, the equality constraints are violated in \texttt{SL+LD} and \texttt{LDF}. This is due to the errors of NN hinder ensuring the equality constraints unless equality embedding is used (see Proposition~\ref{prop:inexact_equality} and Corollary~\ref{cor:positive_error}). By contrast, \texttt{DeepLDE} achieves not only the smallest optimality gap but also provides feasible solutions for both Case 57 and Case 118.
\subsubsection{Computation time}
As shown in Table~\ref{table:ACOPF}, the computation time for the NN-based approaches is significantly shorter than that of \texttt{PYPOWER}. The computation time for the NN-based approaches without equality embedding (\texttt{SL}, \texttt{SL+LD}, and \texttt{LDF}) is shorter than that of \texttt{DeepLDE}. However they cannot satisfy the nonlinear power flow equations (\ref{eq:active_pf})-(\ref{eq:reactive_pf}), and thus the solutions are highly infeasible. Indeed, there is a tradeoff between feasibility and computation time as shown in Fig.~\ref{fig:case118_tradeoff}. By contrast, the solutions of \texttt{DeepLDE} are feasible and about 20 times faster than \texttt{PYPOWER} in Case 57 and Case 118. The computation times for the NN-based methods with equality embedding (\texttt{SL+E}, \texttt{SL+LDE}, and \texttt{DeepLDE}) are all similar except \texttt{DC3}. Furthermore, \texttt{DC3} becomes slower as the size of network grows, implying that \texttt{DC3} is more sensitive to the scale of the problem; this is because the correction method in \texttt{DC3} calculates the Jacobian matrix from entire inequality constraints every correction iteration. By contrast, the primal-dual method in \texttt{DeepLDE} is applied only during the training process. As a result, \texttt{DeepLDE} maintains fast computation even when dealing with large-scale problems.
\begin{figure}[t]
     \centering
     \begin{subfigure}[b]{0.99\columnwidth}
     \centering
         \includegraphics[width=\columnwidth]{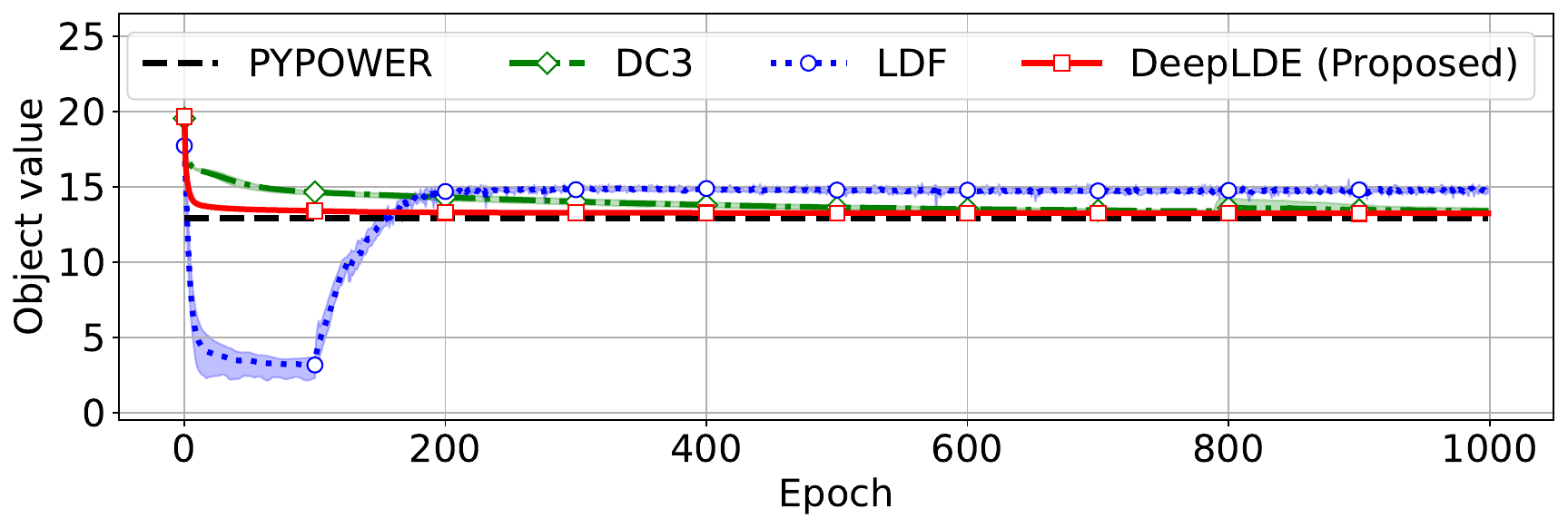}
         \caption{Object funtion}
         \label{fig:obj_function}
     \end{subfigure}
     \begin{subfigure}[b]{0.99\columnwidth}
     \centering
         \includegraphics[width=\columnwidth]{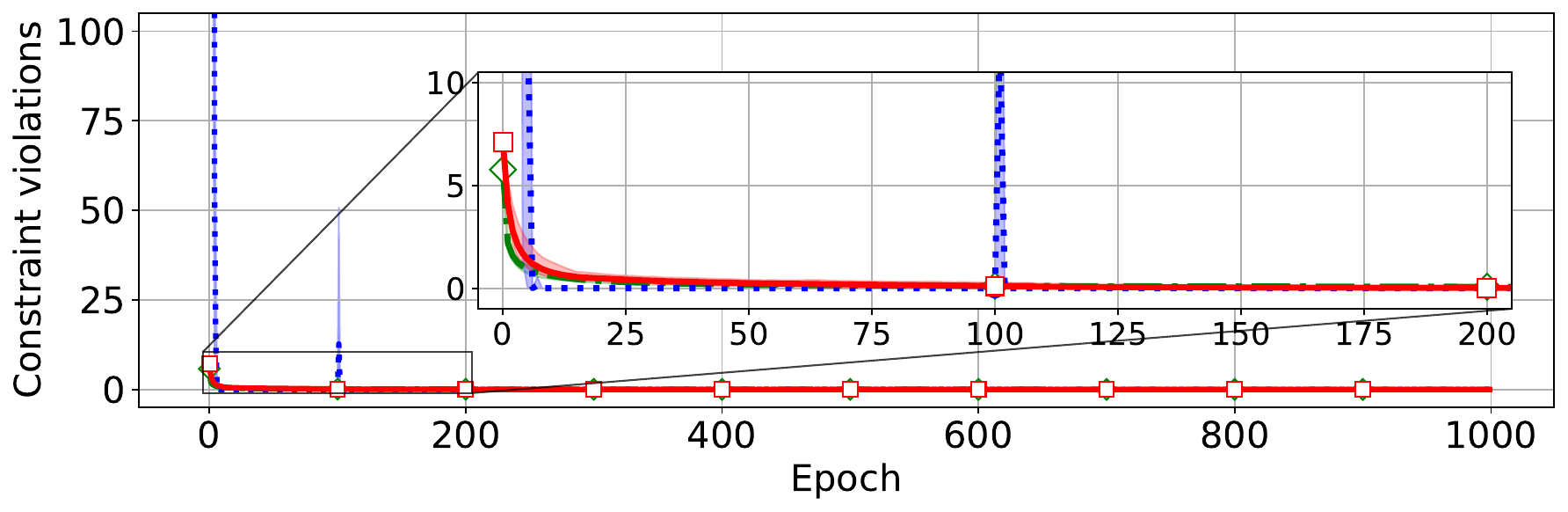}
         \caption{Maximum inequality constraint violations}
         \label{fig:const_violation_ineq}
     \end{subfigure}
     \begin{subfigure}[b]{0.99\columnwidth}
     \centering
         \includegraphics[width=\columnwidth]{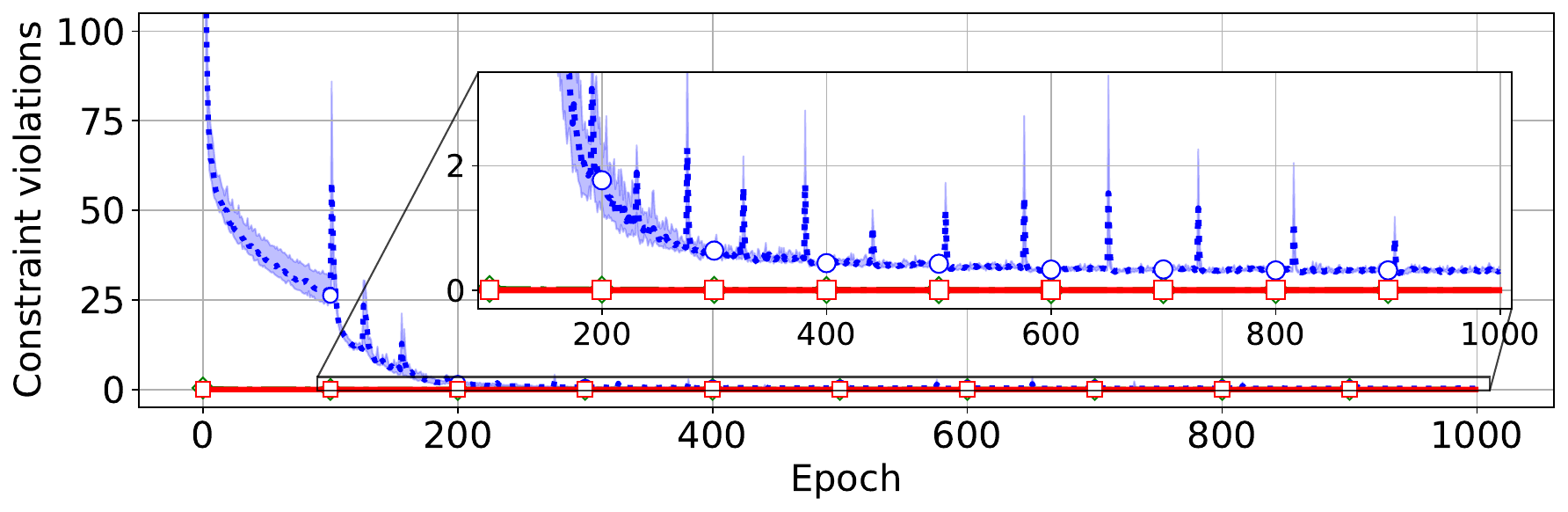}
         \caption{Maximum equality constraint violations}
         \label{fig:const_violation_eq}
     \end{subfigure}
     \caption{Optimality and feasibility of \texttt{DC3}, \texttt{LDF}, and \texttt{DeepLDE}. For \texttt{LDF} and \texttt{DeepLDE}, we use the sum of every inner and warm-up iteration as an epoch as described in Table~\ref{table:iteration}.}
     \label{fig:EE_comparison}
\end{figure}

\subsection{Benefit of Equality Embedding}
\label{sec:exp_eliminiation}
To evaluate the benefit of equality embedding, we compare the learning curves and the maximum constraint violations of \texttt{DC3}, \texttt{LDF}, and \texttt{DeepLDE} using validation data. As illustrated in Fig.~\ref{fig:obj_function}, the objective functions decrease gradually in \texttt{DC3} and \texttt{DeepLDE} as training epochs increase. However, it is not the case of \texttt{LDF}. At the beginning of the training, the object function of \texttt{LDF} decreases even below that of \texttt{PYPOWER}. However, it quickly increases again, and eventually becomes the poorest compared to all other methods. This odd but interesting result is because the errors of \texttt{LDF} hinder ensuring the equality constraints as discussed in Proposition~\ref{prop:inexact_equality} and Corollary~\ref{cor:positive_error}. In contrast to the rapid decrement of inequality constraint violation in all methods, as illustrated in Fig.~\ref{fig:const_violation_ineq}, Fig.~\ref{fig:const_violation_eq} shows that the equality constraints of \texttt{LDF} are never met. Thus, the Lagrange multiplier of the equality constraints $\bm{\mu}$ is continuously updated due to the equality violation. As a result, the weights of \texttt{LDF} are updated in a way that prioritizes reducing equality constraint violations rather than optimizing the object function (\ref{eq:generation_cost}). By contrast, \texttt{DeepLDE} remains feasible while achieving the smallest optimality gap among all other methods in Fig.~\ref{fig:EE_comparison}.

\section{Conclusions}
\label{sec:conclusion}
In this paper, we proposed a framework of deep Lagrange dual with equality embedding (\texttt{DeepLDE}) that successfully integrates equality and inequality constraints within the NNs. To satisfy the equality constraints, we let the NNs infer a part of solution and determine the rest of solution using implicit layers. To satisfy the inequality constraints, we applied the primal-dual method, which does not require any post-processing (unlike \texttt{DC3}). In doing this, we proved that \texttt{DeepLDE} converges. Furthermore, by leveraging the concept of $\varepsilon$-inexact solution, we showed that, if equality embedding is skipped (unlike \texttt{LDF}), primal-dual method never satisfies equality constraints. To verify the performance, we rigorously compared \texttt{DeepLDE} with the conventional solvers and popular other NN-based approaches in solving constrained convex, non-convex, and AC-OPF problems. The results about feasibility, optimality, and computation time confirm that the proposed method always finds a feasible solution that also has the smallest optimality gap among all the NN-based approaches. Furthermore, the computation of \texttt{DeepLDE} is faster than \texttt{DC3} and PYPOWER by 5 and 20 times, respectively. Lastly, \texttt{DeepLDE} aims to address a wide range of real-world, large-scale, constrained optimization problems, which remains as future work.

\appendix
\subsection{Proof of Proposition~\ref{prop:local_convergence}}
\label{subsec:proof_converge}
\begin{proof}
    Let $\rho_t$ be the step size of the Lagrange multiplier $\bm{\lambda}$ at $t^{th}$ outer iteration. Let $r_t$ be the ratio of $\mathbf{w}$ update to $\bm{\lambda}$ update, which is defined as
\begin{equation}
    r_t = \dfrac{\eta \mathcal{N}_\mathbf{w}(t)}{\rho_t},
    \label{eq:proof_ratio}
\end{equation}
where $\mathcal{N}_\mathbf{w}(t)$ denotes the total number of SGD steps to update $\mathbf{w}$ at $t^{th}$ outer iteration of Algorithm~\ref{alg:DeepLDE_noW}. Let $I_t$ denote the maximum inner iteration at $t^{th}$ outer iteration. Thus, $\mathcal{N}_\mathbf{w}(t) = |\mathcal{D}|I_t$ where $|\mathcal{D}|$ is the cardinality of $\mathcal{D}$, and $I_{t+1} = I_t + \beta,\forall t\geq 0$ and thus
\begin{equation}
    I_t = I_0 + \beta t,\quad t\geq 0.
    \label{eq:proof_general}
\end{equation}
Following the theorem from \cite{nandwani2019primal}, Algorithm~\ref{alg:DeepLDE_noW} converges to a local min-max point if $\Lim{t\rightarrow\infty}r_t = \infty$. Since $\rho_t = \frac{\rho_0}{1+\gamma t}$, we have
\begin{equation}
\begin{aligned}
    r_t
    &=\begin{cases}
        \eta|\mathcal{D}|\dfrac{(1+\gamma t)(I_0  + \beta t)}{\rho_0},&\text{if $\beta>0$ and $\gamma>0$};\\
        \eta|\mathcal{D}|\dfrac{I_0  + \beta t}{\rho_0},&\text{if $\beta>0$ and $\gamma=0$};\\
        \eta|\mathcal{D}|\dfrac{(1+\gamma t)I_0 }{\rho_0},&\text{if $\beta=0$ and $\gamma>0$};\\
        \eta|\mathcal{D}|\dfrac{I_0 }{\rho_0},&\text{if $\beta=0$ and $\gamma=0$}.
    \end{cases}
\end{aligned}
\label{eq:proof_rt}
\end{equation}
Note that we do not consider the case when $\beta$ or $\gamma$ is negative. Thus, $\Lim{t\rightarrow\infty}r_t = \infty$ except for the case when both $\beta$ and $\gamma$ are 0, which implies Algorithm~\ref{alg:DeepLDE_noW} converges if $\beta>0$ or $\gamma>0$.
\end{proof}

\subsection{Proof of Proposition~\ref{prop:inexact_equality}}
\label{subsec:inexact_equality}
\begin{proof}
Since $\mathbf{1}^{\intercal}S(h(\mathbf{y}_\varepsilon)) = ||h(\mathbf{y}_\varepsilon)||_1$, and $\mathbf{y}^*$ is an optimal solution of (\ref{eq:prelim_opt}), $h(\mathbf{y}^*) = \mathbf{0}$. Then, we have
\begin{equation}
\begin{aligned}
    \mathbb{E}[\mathbf{1}^{\intercal}S(h(\mathbf{y}_\varepsilon))] &= 
    \mathbb{E}[||h(\mathbf{y}_\varepsilon)||_1]\\ &= \mathbb{E}[||h(\mathbf{y}_\varepsilon) - h(\mathbf{y}^*)||_1]\\
    &= \mathbb{E}[||h(\mathbf{y^*}+\varepsilon) - h(\mathbf{y}^*)||_1].
\end{aligned}
\label{eq:proof_ExpYEps}
\end{equation}
Now, we linearize RHS of (\ref{eq:proof_ExpYEps}) with respect to $\mathbf{y}^*$ for small $\varepsilon$, and obtain
\begin{equation}
\begin{aligned}
    \mathbb{E}[||h(\mathbf{y^*}+\varepsilon) - h(\mathbf{y}^*)||_1] &\simeq \mathbb{E}[||\mathbf{J}_h(\mathbf{y}^*)\varepsilon||_1]\\
    &=\sum_{i=1}^{n_{\text{eq}}}\mathbb{E}[|\{\mathbf{J}_h(\mathbf{y}^*)\varepsilon\}_i|],
\end{aligned}
\label{eq:proof_linearization}
\end{equation}
where $\{\cdot\}_i$ returns $i^{th}$ element of the given vector. Since $\varepsilon\sim\mathcal{N}(\mathbf{0}, \sigma^2\mathbf{I})$ by Assumption~\ref{assume:normal_error}, $\mathbf{J}_h\varepsilon\sim\mathcal{N}(\mathbf{0}, \sigma^2\mathbf{J}_h \mathbf{J}_h^{\intercal})$. Thus, $\{\mathbf{J}_h\varepsilon\}_i\sim\mathcal{N}(\mathbf{0}, \sigma^2\{\mathbf{J}_h \mathbf{J}_h^{\intercal}\}_{ii})$ where $\{\cdot\}_{ii}$ returns $i^{th}$ row and $i^{th}$ column element of the given matrix. Also, RHS of (\ref{eq:proof_linearization}) is
\begin{equation}
\begin{aligned}
    \sum_{i=1}^{n_{\text{eq}}}\mathbb{E}[|\{\mathbf{J}_h(\mathbf{y}^*)\varepsilon\}_i|] 
    &= \sum_{i=1}^{n_{\text{eq}}}\sqrt{\dfrac{2\sigma^2\{\mathbf{J}_h(\mathbf{y}^*)\mathbf{J}_h(\mathbf{y}^*)^{\intercal}\}_{ii}}{\pi}}\\
    &=\sqrt{\dfrac{2}{\pi}}\sigma||\mathbf{J}_h(\mathbf{y}^*)||_*,
\end{aligned}
\end{equation}
where $||\mathbf{J}_h||_* = \tr(\sqrt{\mathbf{J}_h\mathbf{J}_h^{\intercal}})$. Hence, $\mathbb{E}[\mathbf{1}^{\intercal}S(h(\mathbf{y}^*))] = \sqrt{\frac{2}{\pi}}\sigma||\mathbf{J}_h(\mathbf{y}^*)||_*$ holds, and the proof of Proposition~\ref{prop:inexact_equality} is complete.
\end{proof}

\subsection{Proof of Corollary~\ref{cor:positive_error}}
\label{subsec:positive_error}
\begin{proof}
By Proposition~\ref{prop:inexact_equality}, $\mathbb{E}[\mathbf{1}^{\intercal}S(h(\mathbf{y}^*))] = \sqrt{\frac{2}{\pi}}\sigma||\mathbf{J}_h(\mathbf{y}^*)||_*$ where $||\mathbf{J}_h||_* = \tr(\sqrt{\mathbf{J}_h\mathbf{J}_h^{\intercal}})$. Note that $||\mathbf{J}_h||_*$ is the sum of eigenvalues of $\sqrt{\mathbf{J}_h\mathbf{J}_h^{\intercal}}$. Since $\sqrt{\mathbf{J}_h\mathbf{J}_h^{\intercal}}$ is a symmetric positive semidefinite matrix, all the eigenvalues of $\sqrt{\mathbf{J}_h\mathbf{J}_h^{\intercal}}$ are non-negative. Also, $\mathbf{J}_h(\mathbf{y}^*) \neq \mathbf{0}$ implies $\rank(\sqrt{\mathbf{J}_h\mathbf{J}_h^{\intercal}})>0$, which means $\sqrt{\mathbf{J}_h\mathbf{J}_h^{\intercal}}$ has at least one non-zero positive eigenvalue. Hence, $\mathbb{E}[\mathbf{1}^{\intercal}S(h(\mathbf{y}^*))] > 0$ holds, and the proof of Corollary~\ref{cor:positive_error} is complete.
\end{proof}

\subsection{Data generation process}
\label{appendix:data-generate}
We describe the data generation process of linear constrained programming (\ref{eq:QP})$-$(\ref{eq:NonConv}), and AC-OPF (\ref{eq:AC-OPF}). Note that the dataset is divided into training, validation, and test in a 10:1:1 ratio as done in \cite{donti2021DC3}.
\subsubsection{Linear constrained nonlinear programming}
We sample the elements of $\mathbf{A}$ and $\mathbf{G}$ in (\ref{eq:QP})$-$(\ref{eq:NonConv}) from a unit normal distribution. Let $\mathbf{A}^+$ be the Moore-Penrose pseudo-inverse of $\mathbf{A}$. Let $\{\mathbf{d}\}_j$ be $j$th elements of $\mathbf{d}$, and assume that $\{\mathbf{d}\}_j\in[-1,1]$. Since the equality constraints of (\ref{eq:QP})$-$(\ref{eq:NonConv}) are expressed by $\mathbf{A}\mathbf{A}^+\mathbf{d} = \mathbf{d}$, we have
\begin{equation}
    \mathbf{G}\mathbf{A}^+\mathbf{d} \leq \sum_{j}|\{\mathbf{G}\mathbf{A}^+\}_{ij}|
\end{equation}
where $\{\cdot\}_{ij}$ denotes $i^{th}$ row and $j$th column of given matrix \cite{donti2021DC3}. Now, let $\mathbf{h} = \sum_{j}|\{\mathbf{G}\mathbf{A}^+\}_{ij}|$. As a result, the conventional solvers can solve (\ref{eq:QP})$-$(\ref{eq:NonConv}), and ensure feasbile solutions for every $\{\mathbf{d}\}_j\in[-1,1]$. Hence, we generate 10,000 data by sampling $\mathbf{d}$ from a uniform distribution on [-1,1].
\subsubsection{AC-OPF}
We sample $\mathbf{p}_d$ from the following truncated normal distribution \cite{zamzam2020learning, donti2021DC3}:
\begin{equation}
    \mathbf{p}_d \sim \mathcal{TN}(\mathbf{p}_d^0, \bm{\Sigma}_{\mathbf{p}}, (1-\mu)\mathbf{p}_d^0 ,(1+\mu)\mathbf{p}_d^0)
\end{equation}
where $\mathbf{p}_d^0$ denotes a base active power demand provided in MATPOWER \cite{zimmerman2010matpower}, and $\bm{\Sigma}_{\mathbf{p}}$ denotes the covariance matrix which represents relationship between the patterns of active power demand on the buses. We set the correlation parameters of $\bm{\Sigma}_{\mathbf{p}}$ to 0.5 and $\mu$ to 0.7 as done in \cite{donti2021DC3}. In the case of $\mathbf{q}_d$, we sample the power factor of each load from a uniform distribution with interval $[0.8, 1.0]$. Thus, we obtain $\mathbf{q}_d$ from the sampled power factors and the generated $\mathbf{p}_d$. After that, we discard the demand samples that generate infeasible solutions. Hence, we use 2,400 feasible samples of power demand for training, validation, and test.

\bibliographystyle{IEEEtran}
\bibliography{IEEEexample}

\begin{thebibliography}{10}
\providecommand{\url}[1]{#1}
\csname url@samestyle\endcsname
\providecommand{\newblock}{\relax}
\providecommand{\bibinfo}[2]{#2}
\providecommand{\BIBentrySTDinterwordspacing}{\spaceskip=0pt\relax}
\providecommand{\BIBentryALTinterwordstretchfactor}{4}
\providecommand{\BIBentryALTinterwordspacing}{\spaceskip=\fontdimen2\font plus
\BIBentryALTinterwordstretchfactor\fontdimen3\font minus
  \fontdimen4\font\relax}
\providecommand{\BIBforeignlanguage}[2]{{%
\expandafter\ifx\csname l@#1\endcsname\relax
\typeout{** WARNING: IEEEtran.bst: No hyphenation pattern has been}%
\typeout{** loaded for the language `#1'. Using the pattern for}%
\typeout{** the default language instead.}%
\else
\language=\csname l@#1\endcsname
\fi
#2}}
\providecommand{\BIBdecl}{\relax}
\BIBdecl

\bibitem{kim2024unsupervised}
M.~Kim and H.~Kim, ``{Unsupervised Deep Lagrange Dual with Equation Embedding
  for AC Optimal Power Flow},'' \emph{IEEE Transactions on Power Systems},
  2024.

\bibitem{boyd2004convex}
S.~Boyd, S.~P. Boyd, and L.~Vandenberghe, \emph{{Convex Optimization}}.\hskip
  1em plus 0.5em minus 0.4em\relax Cambridge university press, 2004.

\bibitem{hasan2020survey}
F.~Hasan, A.~Kargarian, and A.~Mohammadi, ``{A Survey on Applications of
  Machine Learning For Optimal Power Flow},'' in \emph{2020 IEEE Texas Power
  and Energy Conference (TPEC)}.\hskip 1em plus 0.5em minus 0.4em\relax IEEE,
  2020, pp. 1--6.

\bibitem{bengio2021machine}
Y.~Bengio, A.~Lodi, and A.~Prouvost, ``{Machine Learning for Combinatorial
  Optimization: a Methodological Tour d’Horizon},'' \emph{European Journal of
  Operational Research}, vol. 290, no.~2, pp. 405--421, 2021.

\bibitem{chen2022learning}
Y.~Chen, L.~Zhang, and B.~Zhang, ``{Learning to Solve DCOPF: A Duality
  Approach},'' \emph{Electric Power Systems Research}, vol. 213, p. 108595,
  2022.

\bibitem{pan2020deepopf}
X.~Pan, T.~Zhao, M.~Chen, and S.~Zhang, ``{DeepOPF: A Deep Neural Network
  Approach for Security-constrained DC Optimal Power Flow},'' \emph{IEEE
  Transactions on Power Systems}, vol.~36, no.~3, pp. 1725--1735, 2020.

\bibitem{nellikkath2021physics}
R.~Nellikkath and S.~Chatzivasileiadis, ``{Physics-informed Neural Networks for
  Minimising Worst-case Violations in DC Optimal Power Flow},'' in \emph{2021
  IEEE International Conference on Communications, Control, and Computing
  Technologies for Smart Grids (SmartGridComm)}.\hskip 1em plus 0.5em minus
  0.4em\relax IEEE, 2021, pp. 419--424.

\bibitem{nellikkath2022physics}
------, ``{Physics-informed Neural Networks for AC Optimal Power Flow},''
  \emph{Electric Power Systems Research}, vol. 212, p. 108412, 2022.

\bibitem{nandwani2019primal}
Y.~Nandwani, A.~Pathak, and P.~Singla, ``A {P}rimal {D}ual {F}ormulation for
  {D}eep {L}earning with {C}onstraints,'' \emph{Advances in Neural Information
  Processing Systems}, vol.~32, 2019.

\bibitem{fioretto2020Lagrangian}
F.~Fioretto, P.~V. Hentenryck, T.~W. Mak, C.~Tran, F.~Baldo, and M.~Lombardi,
  ``Lagrangian {D}uality for {C}onstrained {D}eep {L}earning,'' in \emph{Joint
  European Conference on Machine Learning and Knowledge Discovery in
  Databases}.\hskip 1em plus 0.5em minus 0.4em\relax Springer, 2020, pp.
  118--135.

\bibitem{fioretto2020predicting}
F.~Fioretto, T.~W. Mak, and P.~Van~Hentenryck, ``{Predicting AC Optimal Power
  Flows: Combining Deep Learning and Lagrangian Dual Methods},'' in
  \emph{Proceedings of the AAAI Conference on Artificial Intelligence},
  vol.~34, no.~01, 2020, pp. 630--637.

\bibitem{kotary2022fast}
J.~Kotary, F.~Fioretto, and P.~Van~Hentenryck, ``{Fast Approximations for Job
  Shop Scheduling: A Lagrangian Dual Deep Learning Method},'' in
  \emph{Proceedings of the AAAI Conference on Artificial Intelligence},
  vol.~36, no.~7, 2022, pp. 7239--7246.

\bibitem{park2022self}
S.~Park and P.~Van~Hentenryck, ``{Self-Supervised Primal-Dual Learning for
  Constrained Optimization},'' \emph{arXiv preprint arXiv:2208.09046}, 2022.

\bibitem{chen2021enforcing}
B.~Chen, P.~L. Donti, K.~Baker, J.~Z. Kolter, and M.~Berg{\'e}s, ``{Enforcing
  Policy Feasibility Constraints through Differentiable Projection for Energy
  Optimization},'' in \emph{Proceedings of the Twelfth ACM International
  Conference on Future Energy Systems}, 2021, pp. 199--210.

\bibitem{kim2022projection}
M.~Kim and H.~Kim, ``{Projection-aware Deep Neural Network for DC Optimal Power
  Flow Without Constraint Violations},'' in \emph{2022 IEEE International
  Conference on Communications, Control, and Computing Technologies for Smart
  Grids (SmartGridComm)}.\hskip 1em plus 0.5em minus 0.4em\relax IEEE, 2022,
  pp. 116--121.

\bibitem{amos2017optnet}
B.~Amos and J.~Z. Kolter, ``Opt{Net}: {D}ifferentiable {O}ptimization as a
  {L}ayer in {N}eural {N}etworks,'' in \emph{International Conference on
  Machine Learning}.\hskip 1em plus 0.5em minus 0.4em\relax PMLR, 2017, pp.
  136--145.

\bibitem{agrawal2019differentiable}
A.~Agrawal, B.~Amos, S.~Barratt, S.~Boyd, S.~Diamond, and J.~Z. Kolter,
  ``{Differentiable Convex Optimization Layers},'' \emph{Advances in Neural
  Information Processing Systems}, vol.~32, 2019.

\bibitem{donti2021DC3}
P.~Donti, D.~Rolnick, and J.~Z. Kolter, ``{DC3}: A {L}earning {M}ethod for
  {O}ptimization with {H}ard {C}onstraints,'' in \emph{International Conference
  on Learning Representations}, 2021.

\bibitem{babaeinejadsarookolaee2019power}
S.~Babaeinejadsarookolaee, A.~Birchfield, R.~D. Christie, C.~Coffrin,
  C.~DeMarco, R.~Diao, M.~Ferris, S.~Fliscounakis, S.~Greene, R.~Huang
  \emph{et~al.}, ``{The Power Grid Library for Benchmarking AC Optimal Power
  Flow Algorithms},'' \emph{arXiv preprint arXiv:1908.02788}, 2019.

\bibitem{bergstra2011algorithms}
J.~Bergstra, R.~Bardenet, Y.~Bengio, and B.~K{\'e}gl, ``{Algorithms for
  Hyper-parameter Optimization},'' \emph{Advances in Neural Information
  Processing Systems}, vol.~24, 2011.

\bibitem{hutter2011sequential}
F.~Hutter, H.~H. Hoos, and K.~Leyton-Brown, ``{Sequential Model-based
  Optimization for General Algorithm Configuration},'' in \emph{International
  conference on learning and intelligent optimization}.\hskip 1em plus 0.5em
  minus 0.4em\relax Springer, 2011, pp. 507--523.

\bibitem{snoek2012practical}
J.~Snoek, H.~Larochelle, and R.~P. Adams, ``{Practical Bayesian Optimization of
  Machine Learning Algorithms},'' \emph{Advances in Neural Information
  Processing Systems}, vol.~25, 2012.

\bibitem{bergstra2012random}
J.~Bergstra and Y.~Bengio, ``{Random Search for Hyper-parameter
  Optimization.}'' \emph{Journal of machine learning research}, vol.~13, no.~2,
  2012.

\bibitem{li2017learning}
K.~Li and J.~Malik, ``{Learning to Optimize},'' in \emph{International
  Conference on Learning Representations}, 2017.

\bibitem{chenL2O2022learning}
T.~Chen, X.~Chen, W.~Chen, Z.~Wang, H.~Heaton, J.~Liu, and W.~Yin, ``{Learning
  to Optimize: A Primer and a Benchmark},'' \emph{The Journal of Machine
  Learning Research}, vol.~23, no.~1, pp. 8562--8620, 2022.

\bibitem{donti2017task}
P.~Donti, B.~Amos, and J.~Z. Kolter, ``{Task-based End-to-End Model Learning in
  Stochastic Optimization},'' \emph{Advances in neural information processing
  systems}, vol.~30, 2017.

\bibitem{wilder2019end}
B.~Wilder, E.~Ewing, B.~Dilkina, and M.~Tambe, ``{End to End Learning and
  Optimization on Graphs},'' \emph{Advances in Neural Information Processing
  Systems}, vol.~32, 2019.

\bibitem{poganvcic2019differentiation}
M.~V. Pogan{\v{c}}i{\'c}, A.~Paulus, V.~Musil, G.~Martius, and M.~Rolinek,
  ``{Differentiation of Blackbox Combinatorial Solvers},'' in
  \emph{International Conference on Learning Representations}, 2019.

\bibitem{elmachtoub2022smart}
A.~N. Elmachtoub and P.~Grigas, ``{Smart “Predict, then Optimize”},''
  \emph{Management Science}, vol.~68, no.~1, pp. 9--26, 2022.

\bibitem{kotary2021end}
J.~Kotary, F.~Fioretto, P.~Van~Hentenryck, and B.~Wilder, ``{End-to-end
  Constrained Optimization Learning: A Survey},'' \emph{arXiv preprint
  arXiv:2103.16378}, 2021.

\bibitem{sun2018learning}
H.~Sun, X.~Chen, Q.~Shi, M.~Hong, X.~Fu, and N.~D. Sidiropoulos, ``{Learning to
  Optimize: Training Deep Neural Networks for Interference Management},''
  \emph{IEEE Transactions on Signal Processing}, vol.~66, no.~20, pp.
  5438--5453, 2018.

\bibitem{zamzam2020learning}
A.~S. Zamzam and K.~Baker, ``Learning {O}ptimal {S}olutions for {E}xtremely
  {F}ast {AC} {O}ptimal {P}ower {F}low,'' in \emph{2020 IEEE International
  Conference on Communications, Control, and Computing Technologies for Smart
  Grids (SmartGridComm)}.\hskip 1em plus 0.5em minus 0.4em\relax IEEE, 2020,
  pp. 1--6.

\bibitem{zhang2021convex}
L.~Zhang, Y.~Chen, and B.~Zhang, ``{A Convex Neural Network Solver for DCOPF
  with Generalization Guarantees},'' \emph{IEEE Transactions on Control of
  Network Systems}, 2021.

\bibitem{venzke2020learning}
A.~Venzke, G.~Qu, S.~Low, and S.~Chatzivasileiadis, ``{Learning Optimal Power
  Flow: Worst-case Guarantees for Neural Networks},'' in \emph{2020 IEEE
  International Conference on Communications, Control, and Computing
  Technologies for Smart Grids (SmartGridComm)}.\hskip 1em plus 0.5em minus
  0.4em\relax IEEE, 2020, pp. 1--7.

\bibitem{ul2022learning}
Z.~ul~Abdeen, H.~Yin, V.~Kekatos, and M.~Jin, ``{Learning Neural Networks under
  Input-Output Specifications},'' in \emph{2022 American Control Conference
  (ACC)}.\hskip 1em plus 0.5em minus 0.4em\relax IEEE, 2022, pp. 1515--1520.

\bibitem{zhao2023ensuring}
T.~Zhao, X.~Pan, M.~Chen, and S.~Low, ``{Ensuring {DNN} Solution Feasibility
  for Optimization Problems with Linear Constraints},'' in \emph{The Eleventh
  International Conference on Learning Representations}, 2023.

\bibitem{ling2020solving}
Z.~Ling, X.~Tao, Y.~Zhang, and X.~Chen, ``{Solving Optimization Problems
  Through Fully Convolutional Networks: An Application to the Traveling
  Salesman Problem},'' \emph{IEEE Transactions on Systems, Man, and
  Cybernetics: Systems}, vol.~51, no.~12, pp. 7475--7485, 2020.

\bibitem{ling2021can}
Z.~Ling, R.~Liu, Y.~Zhang, and X.~Chen, ``{Can Deep Learning Solve Parametric
  Mathematical Programming? An Application to 0--1 Linear Programming Through
  Image Representation},'' \emph{IEEE Transactions on Systems, Man, and
  Cybernetics: Systems}, vol.~52, no.~9, pp. 5656--5667, 2021.

\bibitem{tabas2022computationally}
D.~Tabas and B.~Zhang, ``Computationally {E}fficient {S}afe {R}einforcement
  {L}earning for {P}ower {S}ystems,'' in \emph{2022 American Control Conference
  (ACC)}.\hskip 1em plus 0.5em minus 0.4em\relax IEEE, 2022, pp. 3303--3310.

\bibitem{li2022learning}
M.~Li, S.~Kolouri, and J.~Mohammadi, ``{Learning to Solve Optimization Problems
  with Hard Linear Constraints},'' \emph{arXiv preprint arXiv:2208.10611},
  2022.

\bibitem{blanchini2008set}
F.~Blanchini and S.~Miani, \emph{{Set-theoretic Methods in Control}}.\hskip 1em
  plus 0.5em minus 0.4em\relax Springer, 2008, vol.~78.

\bibitem{he2015neural}
W.~He, A.~O. David, Z.~Yin, and C.~Sun, ``{Neural Network Control of a Robotic
  Manipulator with Input Deadzone and Output Constraint},'' \emph{IEEE
  Transactions on Systems, Man, and Cybernetics: Systems}, vol.~46, no.~6, pp.
  759--770, 2015.

\bibitem{fontaine2014constraint}
D.~Fontaine, LaurentMichel, and P.~Van~Hentenryck, ``{Constraint-based
  Lagrangian Relaxation},'' in \emph{Principles and Practice of Constraint
  Programming: 20th International Conference, CP 2014, Lyon, France, September
  8-12, 2014. Proceedings 20}.\hskip 1em plus 0.5em minus 0.4em\relax Springer,
  2014, pp. 324--339.

\bibitem{nocedal2006numerical}
J.~Nocedal and S.~J. Wright, \emph{Numerical {O}ptimization (2nd
  edition)}.\hskip 1em plus 0.5em minus 0.4em\relax Springer, 2006.

\bibitem{jin2020local}
C.~Jin, P.~Netrapalli, and M.~Jordan, ``What is {L}ocal {O}ptimality in
  {N}onconvex-nonconcave {M}inimax {O}ptimization?'' in \emph{International
  Conference on Machine Learning}.\hskip 1em plus 0.5em minus 0.4em\relax PMLR,
  2020, pp. 4880--4889.

\bibitem{stellato2020OSQP}
B.~Stellato, G.~Banjac, P.~Goulart, A.~Bemporad, and S.~Boyd, ``{OSQP}: {A}n
  {O}perator {S}plitting {S}olver for {Q}uadratic {P}rograms,''
  \emph{Mathematical Programming Computation}, vol.~12, no.~4, pp. 637--672,
  2020.

\bibitem{wachter2006implementation}
A.~W{\"a}chter and L.~T. Biegler, ``On the {I}mplementation of an
  {I}nterior-point {F}ilter {L}ine-search {A}lgorithm for {L}arge-scale
  {N}onlinear {P}rogramming,'' \emph{Mathematical programming}, vol. 106,
  no.~1, pp. 25--57, 2006.

\bibitem{zimmerman2010matpower}
R.~D. Zimmerman, C.~E. Murillo-S{\'a}nchez, and R.~J. Thomas, ``{MATPOWER:
  Steady-state Operations, Planning, and Analysis Tools for Power Systems
  Research and Education},'' \emph{IEEE Transactions on power systems},
  vol.~26, no.~1, pp. 12--19, 2010.

\bibitem{clevert2015fast}
D.-A. Clevert, T.~Unterthiner, and S.~Hochreiter, ``Fast and {A}ccurate {D}eep
  {N}etwork {L}earning by {E}xponential {L}inear {U}nits {(ELUs)},''
  \emph{arXiv preprint arXiv:1511.07289}, 2015.

\bibitem{srivastava2014dropout}
N.~Srivastava, G.~Hinton, A.~Krizhevsky, I.~Sutskever, and R.~Salakhutdinov,
  ``{Dropout: A Simple Way to Prevent Neural Networks from Overfitting},''
  \emph{The journal of machine learning research}, vol.~15, no.~1, pp.
  1929--1958, 2014.

\bibitem{kingma2014adam}
D.~P. Kingma and J.~Ba, ``Adam: {A} {M}ethod for {S}tochastic {O}ptimization,''
  in \emph{International Conference on Learning Representations}, 2015.

\bibitem{pan2022deepopf}
X.~Pan, M.~Chen, T.~Zhao, and S.~H. Low, ``{DeepOPF: A Feasibility-optimized
  Deep Neural Network Approach for AC Optimal Power Flow Problems},''
  \emph{IEEE Systems Journal}, 2022.

\end{thebibliography}

\newpage

\vfill

\end{document}